\documentclass[12pt,a4paper]{amsart}
\usepackage[top=35mm, bottom=35mm, left=30mm, right=30mm]{geometry}
\usepackage{mathptmx}
\usepackage{mathrsfs}
\usepackage{verbatim}
\usepackage{url}
\usepackage[all]{xy}
\usepackage{xcolor}
\usepackage[colorlinks=true,citecolor=blue]{hyperref}
\usepackage{amsmath,amsthm}
\usepackage{amssymb}
\usepackage{ulem}
\usepackage{tikz}
\usetikzlibrary{arrows}
\usetikzlibrary{patterns}
\usetikzlibrary{lindenmayersystems,backgrounds}
\usepackage[colorinlistoftodos,prependcaption,textsize=tiny, textwidth=20mm,]{todonotes}

\makeatletter
\@namedef{subjclassname@2020}{\textup{2020} Mathematics Subject Classification}
\makeatother

\theoremstyle{plain}

\newtheorem{thm}{Theorem}[section]
\newtheorem{thma}{Theorem}

\newtheorem{lem}[thm]{Lemma}
\newtheorem{prop}[thm]{Proposition}
\newtheorem{cor}[thm]{Corollary}
\newtheorem{prob}[thm]{Problem}
\theoremstyle{definition}
\newtheorem{defn}[thm]{Definition}
\newtheorem{rem}[thm]{Remark}

\newcommand{\Z}{{\mathbb{Z}}}
\newcommand{\N}{\mathbb{N}}

\begin{document}
\title{Disjointness with all minimal systems under group actions}

\author[H.~Xu]{Hui Xu}
\address[H. Xu]{Department of Mathematics, Shanghai Normal University, Shanghai 200234, China}
\email{huixu2734@ustc.edu.cn}

\author[X.~Ye]{Xiangdong Ye}
\address[X. Ye]{CAS Wu Wen-Tsun Key Laboratory of Mathematics, University of Science and Technology of China, Hefei, Anhui 230026, China}
\email{yexd@ustc.edu.cn}



\thanks{H. Xu is supported by NNSF 12201599, 12371196 and X. Ye is supported by NNSF 12031019}

\maketitle

\begin{abstract}
Let $G$ be a countable discrete group. We give a necessary and sufficient condition for a transitive $G$-system to be disjoint with all minimal $G$-systems, which implies that if a transitive $G$-system is disjoint with all minimal $G$-systems, then it is $\infty$-transitive, i.e. $(X^k,G)$ is transitive for all $k\in\N$, and has dense minimal points. In addition, we show that any $\infty$-transitive $G$-system with dense distal points are disjoint with all minimal $G$-systems.
\end{abstract}


\section{Introduction}

Let $G$ be a discrete group, the identity element of which is denoted by $e_{G}$. By a $G$-{\it system} (or {\it system} when the group is clear in the context) we mean a compact metrizable space $X$ equipped with a continuous left action of $G$ by homeomorphisms. Such a $G$-system is denoted by $(X, G)$ or $G\curvearrowright X$.  Two $G$-systems $(X,G)$ and $(Y,G)$ are said to be {\it disjoint} if $X\times Y$ is the only nonempty closed subset of $X\times Y$ that is invariant under the diagonal action of $G$ and has full projection onto each coordinate. Such notion of disjointness was introduced by Furstenberg in his  seminal work \cite{Fur}.

 \medskip
It is easy to check that if two systems are disjoint then one of them is minimal (see \cite[Theorem II.1]{Fur}).  Further, Furstenberg showed in \cite{Fur} that every totally transitive $\mathbb{Z}$-system with dense periodic points is disjoint with all minimal $\mathbb{Z}$-systems and every weakly mixing $\mathbb{Z}$-system is disjoint with all minimal distal $\mathbb{Z}$-systems. In particular, the Bernoulli shift on $\{0,1\}^{\mathbb{Z}}$ is disjoint with all minimal $\mathbb{Z}$-systems. Then the  following natural question is left.
\begin{prob}\cite[Problem G]{Fur}\label{Fur's problem}
Describe the $\mathbb{Z}$-system that is disjoint with all minimal $\mathbb{Z}$-systems.
\end{prob}

In \cite{HY}, Huang and Ye showed that a transitive $\mathbb{Z}$-system that is disjoint with all minimal $\mathbb{Z}$-systems has to be weakly mixing and has a dense set of minimal points, and every weakly mixing $\mathbb{Z}$-system with dense small periodic sets or dense regular minimal points is disjoint will all minimal $\mathbb{Z}$-systems. In addition, using the notion of the $m$-set,  they also gave a necessary and sufficient condition for a transitive $\mathbb{Z}$-system to be disjoint with all minimal $\mathbb{Z}$-systems.  Later, it was showed in \cite{DSY, Opr10} that every weakly mixing $\mathbb{Z}$-system with dense distal points is disjoint will all minimal $\mathbb{Z}$-systems, and the relation between $(X,T)$
and its induced spaces for disjointness was discussed in \cite{LYY}.

\medskip
Oprocha gave in \cite{Opr19} a more effective sufficient condition that guarantees a transitive $\mathbb{Z}$-system to be disjoint with all minimal $\mathbb{Z}$-systems. Later, Huang, Shao and Ye gave a complete answer to  Furstenberg's problem   in \cite{HSY} by establishing the necessity of Oprocha's condition, which was left as a question in \cite{Opr19}. However, there are seldom results on  Problem \ref{Fur's problem} if we consider  general group actions.  For a transitive $\mathbb{Z}^{d}$-system, Yu showed in \cite{Yu} that it has a dense set of minimal points if it is disjoint with all minimal $\mathbb{Z}^{d}$-systems. In addition, he also showed a weakly mixing $G$-system with dense set of distal points is disjoint with all minimal $G$-systems if $G$ is an abelian group.

\medskip

The first remarkable result on general group actions was given by Glasner, Tsankov, Weiss and Zucker in \cite{GTWZ}, where the authors showed that for any infinite discrete group $G$, the Bernoulli shift $2^{G}$ is disjoint with all minimal $G$-systems (see \cite{Ber} for a short proof).
The main goal of this paper is to study Problem \ref{Fur's problem} for general countable discrete group actions. We will show that many results holding for $\mathbb{Z}$-actions are also true for   countable discrete group actions. We say a $G$-system $(X,G)$ is $k$-transitive with $k\geq 1$ if the diagonal action of $G$ on the product space $X^{k}=X\times\cdots\times X$ with $k$ copies of $X$ is 
transitive. Further, we say that $(X,G)$ is $\infty$-transitive if it is  $k$-transitive for any $k\geq 1$.

\medskip
The followings are our main results.

\begin{thma}\label{thmA}
Let $G$ be a countable discrete group and  $(X,G)$ be a transitive $G$-system. Then  $(X,G)$ is disjoint with all minimal $G$-systems if and only if
$(X, G)$ is   $\infty$-transitive and there is a countable dense subset $D\subset X$ consisting of minimal points such that for every minimal $G$-system $(Y,G)$,  any $y\in Y$ and any open neighborhood $V$ of $y$, and any nonempty open set $U\subset X$ , there is some $x\in D\cap U$ such that   $N_{G}((x,y), U\times V)=\{g\in G: gx\in U, gy\in V\}$ is right syndetic.

Particularly, if a transitive $G$-system  is disjoint with all minimal $G$-systems, then it is  $\infty$-transitive with dense minimal points.
\end{thma}

The proof of Theorem \ref{thmA} heavily relies on the ideas in \cite{HY, HSY}. In particular,  the authors in \cite{HY} introduced the notion of $m$-set to characterize of disjointness via the recurrence time. In \cite{GTWZ}, the authors introduce the notion of dense orbit set which also characterizes the disjointness and a question is left there that how to characterize the structure of the dense orbit sets in a countable group. This was answered in \cite{KRS} saying that the family of dense orbit sets is the dual to the collection of symmetrically syndetic sets. In Section \ref{char of dis}, we will give some equivalent characterizations of disjointness and study the relations among $m$-sets, dense orbit sets, symmetrically syndetic sets and the minimal points in the Bernoulli shift. We remark that in this paper we only focus on transitive systems disjoint from all minimal systems, see \cite{HY, HSY} for discussions when $(X,G)$ is not transitive for $G=\Z$.

\medskip
Although, the main idea to prove Theorem \ref{thmA} follows from \cite{HY, HSY}, there are lots of essential difficulties needed to overcome for general group actions, and we will show more so that it is more adaptable for applications (see Theorem \ref{stronger than thmA}). The most difficult part is to generalize \cite[Theorem 2.4]{HY} that is to show the existence of $m$-sets in a thickly syndetic set for a countable group. Fortunately, this can be done by strengthening the powerful tools created by Gao-Jackson-Seward in \cite[Chapter 5]{GJS}, which provides some kind of  tilings for any countable groups to some extent. This helps us establish the following result that may have its own interest.
\begin{thma}\label{thmB}
Let $G$ be a countable group and $T$ be a right thickly syndetic subset of $G$. Then there is a nonempty subset $S\subset T$ such that the characteristic function ${\bf 1}_{S}$ of $S$ is a minimal 
point in the Bernoulli shift $2^{G}$.
\end{thma}

In addition, we also have the following sufficient conditions for a transitive system to be disjoint with all minimal systems, which are all known for $\mathbb{Z}$-actions.
\begin{thma}\label{thmC}
Let $G$ be a countable discrete group and  $(X,G)$ be a transitive $G$-system. Then  $(X,G)$ is disjoint with all minimal $G$-systems if it satisfies one of the following conditions:
\begin{itemize}
\item[(1)] $(X,G)$ is totally transitive and has dense small periodic sets;
\item[(2)] $(X,G)$ is totally transitive and has dense regular minimal points;
\item[(3)] $(X,G)$ is $\infty$-transitive and has dense distal points.
\end{itemize}
\end{thma}

  Although Theorem \ref{thmA} provides a necessary and sufficient condition on disjointness from all minimal systems, it is not easy to construct the countable dense set $D$ for a given $G$-system $(X,G)$ generally.  For example, the construction of such $D$ in the Bernoulli shift $2^{G}$ depends on the separated covering property for minimal systems, which is also the main ingredient and difficulty in \cite{GTWZ} to show the Bernoulli disjointness
(for the construction, see Corollary \ref{Bernoulli disjoint}).  For some special countable infinite groups, e.g., residually finite groups, there are dense sets of periodic points in the Bernoulli shift and then the disjointness follows from Theorem \ref{thmC}. Indeed, for a residually finite group $G$ and a finite subset $F\subset G$, there is always a normal subgroup $H$ of $G$ with finite index such that for any $f_1\neq f_2\in F$ lie in distinct cosets of $H$. Then any partial function $\varphi: F\rightarrow\{0,1\}$ can be extended to a periodic point $x$ in $2^{G}$ by setting $x(fh)=\varphi(f)$ for any $f\in F, h\in H$ and $x(g)=0$ otherwise. In addition, it is also necessary to be residually finite for having dense periodic points.

Further, it is not the case that the Bernoulli shift of every countable group admits dense distal points, e.g, Tarski monster groups
(for a proof see Corollary \ref{tarskinotdense}). However, we do not know for which group $G$ the set of distal points is dense in $2^{G}$ beside residually finite groups.

\medskip
Finally, we end this introduction with the following question,  the first item of which is even open for $\mathbb{Z}$-actions (see \cite{DSY, HSY}).
\newtheorem*{theorem}{Question}
\begin{theorem}
\begin{itemize}
\item[(1)] Let $(X,G)$ and $(Y,G)$ be two transitive $G$-systems that are disjoint with all minimal $G$-systems. Is the product system $(X\times Y, G)$ also disjoint with all minimal $G$-systems?
\item[(2)]  For which countable infinite group $G$ is the set of distal points dense in $2^{G}$?
\end{itemize}
\end{theorem}

\medskip
\noindent{\bf Organization of the paper}. In Section 2, we give some basic definitions and lemmas. In Section 3, we give some equivalent characterizations of transitive $G$-systems to be disjoint with all minimal $G$-systems using $m$-sets and dense orbit sets. In Section 4, we prove Theorem \ref{thmA} assuming Proposition \ref{m-set}.
In Section 5, we prove Theorem \ref{thmC}. Finally, we show Proposition \ref{m-set} and Theorem \ref{thmB}.

\section{Preliminary}

In this section we provide necessary notions and results needed in the sequel.

\subsection{Subsets of groups}

\begin{defn}
Let $G$ be a group and $S$ be a subset of $G$. We say that
\begin{itemize}
\item[(1)] $S$ is {\it left syndetic}  (resp. {\it right syndetic}) if there is a finite subset $F\subset G$ with $G=SF$ (resp. $G=FS$);
\item[(2)]  $S$ is {\it left thick}  (resp. {\it right thick}) if for any finite subset $A\subset G$, there is  $g\in G$ with $gA\subset S$ (resp. $Ag\subset S$);
\item[(3)]  $S$ is {\it left piecewise syndetic} (resp. {\it right piecewise syndetic}) if there is a finite subset $F\subset G$ such that for any finite subset $A\subset G$ there is $g\in G$ with $gA\subset SF$ (resp. $Ag\subset FS$);
\item[(4)] $S$ is {\it left thickly syndetic}  (resp. {\it right thickly syndetic}) if for any finite subset $A\subset G$, there is a left syndetic (resp. right syndetic) set $Q_{A}\subset G$ with $Q_AA\subset S$ (resp. $AQ_A\subset S$).
\end{itemize}
\end{defn}
The following lemma is clear from the definitions.
\begin{lem}\label{piecewise synd}
Let $G$ be a group and $S\subset G$. Then $S$ is right piecewise syndetic if and only if $S$ has nonempty intersection with every right thickly syndetic subset of $G$.
\end{lem}
The following are some basic facts on groups and can be found in many textbooks (e.g. \cite[Theorem 3.1]{Roman}).
\begin{lem}\label{group}
Let $G$ be a group and $H$ be a finite index subgroup of $G$. Then
\begin{itemize}
\item[(1)] Every finite index subgroup of $H$ has finite index in $G$;
\item[(2)] If $K$ is also a finite index subgroup of $G$, then $H\cap K$ has finite index in $G$;
\item[(3)] There is a  subgroup $N$ of $H$ that is normal in $G$ and has finite index in $G$.
\end{itemize}
\end{lem}

\begin{lem}\label{syn+thick}
Let $G$ be a group and $S\subset G$ be right thick in $G$. Then
\begin{itemize}
\item[(1)] For any finite index subgroup $H $ of $G$, $S\cap H$ is right thick in $H$;
\item[(2)] For any $g\in G$ and  any finite index subgroup $H $ of $G$, $gSg^{-1}\cap H$ is right thick in $H$.
\end{itemize}
\end{lem}
\begin{proof}
(1) Let $H$ be a finite index subgroup of $G$. Let $G=g_1H\cup g_2H\cup\cdots\cup g_{k}H$ be a cosets decomposition of $G$ with respect to $H$, where $k=[G: H]$ is the index of $H$ in $G$.

In order to show $S\cap H$ is right thick in $H$, we need to show that for any finite subset $A\subset H$, there is some $h\in H$ such that $Ah\subset S\cap H$. Now fix a finite subset $A\subset H$. Let $\widetilde{A}=Ag_1^{-1}\cup \cdots\cup Ag_{k}^{-1}$. Since $S$ is right thick in $G$, there is some $g\in G$ such that $\widetilde{A}g\subset S$. Write $g=g_{i}h$ with $i\in\{1,\ldots, k\}$ and $h\in H$. Then we have
\[Ah=Ag_i^{-1}g_ih\subset \widetilde{A}g\subset S.\]
Clearly, $Ah\subset H$, since $A\subset H$ and $h\in H$. Thus $Ah\subset S\cap H$. This shows that $S\cap H$ is right thick in $H$.

\medskip
(2) According to the clause (1), it suffices to show that $gSg^{-1}$ is also right thick in $G$. While this holds trivially.
\end{proof}

A {\it graph} $\Gamma$ is an ordered pair $(V,E)$ with $V$ being a set and $E$ being a subset of ${V\choose 2}=\{B\subset V: |B|=2\}$. Then $V$ is called the vertex set of $\Gamma$ and $E$ is the edge set of $\Gamma$. For $u,v\in V$, if $\{u,v\}\in E$ then we say that $u,v$ are adjacent and $v$ is a {\it neighbor} of $u$. A subset $I\subset V$ is {\it independent} if no two elements of $I$ are adjacent.  If the vertex set is finite then we call the graph is finite, otherwise it is an infinite graph. The following lemma is well known, e.g. \cite[V. 1]{Bollobas}.
\begin{lem}\label{independent}
Let $\Gamma=(V,E)$ be a finite graph. If every vertex has at most $d$ neighbors, then there is an independent set $I\subset V$ with $|I|\geq |V|/(d+1)$.
\end{lem}

\begin{lem}\label{dis trans}
Let $G$ be a group and $A, B\subset G$. For any $n\in\mathbb{N}$, if $|B|\geq n(|A|^2+1)$, then there are $b_1,\ldots, b_n\in B$ such that $Ab_1,\ldots,Ab_n$ are mutually disjoint.
\end{lem}
\begin{proof}
Define a graph $\Gamma=(V,E)$ with vertex set $V=B$ and two distinct $b, b'\in B$ are adjacent if and only if $Ab\cap Ab'\neq \emptyset$. Thus $b, b'\in B$ are adjacent if and only if $b'\in (A^{-1}Ab\cap B)\setminus\{b\}$. So the number of neighbors of $b$ is no greater than $  |A^{-1}Ab|\leq |A|^{2}$. By Lemma \ref{independent}, there is an independent subset $I\subset B$ with $|I|\geq |B|/(|A|^2+1)\geq n$. That is to say $Ab\cap Ab'=\emptyset$ for any distinct $b,b'\in I$. So we complete the proof.
\end{proof}

We also need the following combinatorial lemma.

\begin{lem}\label{separated syndetic}
Let $G$ be an infinite discrete group and $A\subset G$ be a right syndetic subset. Then for any finite subset $F\subset G$, there is a  right syndetic subset $B$ contained in $A$ such that for any distinct $b, b'\in B\cup\{e_{G}\}$, $Fb\cap Fb'=\emptyset$.
\end{lem}
\begin{proof}
Let $B$ a maximal subset of $A$ under inclusion such that for any distinct $b, b'\in B\cup\{e_{G}\}$, $Fb\cap Fb'=\emptyset$. Then we claim that $A\subset (F^{-1}FB)\cup (F^{-1}F)$. Indeed, it follows from  the maximality of $B$ that for any $a\in A$, there is some $b\in B\cup\{e_{G}\}$ with $Fa\cap Fb\neq \emptyset$ and hence $a\in F^{-1}Fb$. Thus we have $A\subset  (F^{-1}FB)\cup (F^{-1}F)$. Now there is some finite subset $E\subset G$ with $G=EA$, since $A$ is right syndetic. It follows that $$G\subset(EF^{-1}FB)\cup (EF^{-1}F) \subset EF^{-1}F\{e_{G}, b^{-1}\})B,$$
where $b$ is any element in $B$.  Hence $B$ is right syndetic.
\end{proof}

\begin{lem}\label{thick translate}
Let $G$ be a countable discrete group and $T\subset G$ be a right thick subset. Then for any finite subset $F\subset G$, the set $\{g\in G: Fg\subset T\}$ is right thick.
\end{lem}
\begin{proof}
Let $E$ be any finite subset of $G$. Then there is some $h\in G$ such that $FEh\subset T$, since $FE$ is finite. Thus $Eh\subset \{g\in G: Fg\subset T\}$. This show that the set $\{g\in G: Fg\subset T\}$ is right thick.
\end{proof}


\subsection{Topological dynamical system}
Let $G$ be a discrete group. A $G$-{\it system} is a compact metrizable space $X$ equipped a continuous left action of $G$ on $X$, i.e. there is a group homomorphism $\phi$ from $G$ to the homeomorphism group ${\rm Homeo}(X)$ of $X$. More generally, $X$ can be a compact Hausdorff space. But in this paper, we always consider compact metrizable spaces except explicitly mentioned. For brevity, we denote $\phi(g)(x)$ by $gx$ or $g(x)$ or $g\cdot x$, for any $g\in G$ and $x\in X$ and denote the $G$-system by $(X,G)$  or $G\curvearrowright X$. For a subset $Z\subset X$ and a subset $F\subset G$, we denote the orbit of $Z$ under $F$ by $F\cdot Z=\{f(z): f\in F,\ z\in Z\}$ or $FZ$ for short. When $Z=\{x\}$, we denote $F\{x\}$ by $F\cdot x$ or $Fx$ for short.  A subset $Z\subset X$ is $G$-{\it invariant} (or {\it invariant} for short) if $GZ=Z$.  For a closed invariant subset $Z\subset X$, there is an action of $G$ on $Z$ by restrictions and we say $(Z, G)$ is a {\it subsystem} of $(X, G)$. A point $x\in X$ is {\it periodic} if the orbit $Gx$ of $x$ is finite, and is {\it transitive } if $Gx$ is dense in $X$.

\medskip

A $G$-system $(X,G)$ is {\it topologically transitive} (or {\it transitive}) if for any nonempty open subsets $U, V$ of $X$, there is some $g$ such that $gU\cap V\neq\emptyset$. In our setting, $(X,G)$ is transitive if and only if it admits a transitive point.
A $G$-system $(X,G)$ is {\it minimal} if for every point $x\in X$,  the orbit $Gx$ of $x$ under $G$ is dense in $X$.  A point $x\in X$ is a {\it minimal point} if the subsystem $(\overline{Gx}, G)$ is minimal. A point $x\in X$ is an {\it almost periodic point} if for any open neighborhood $U$, the set $\{g\in G: gx\in U\}$ is right syndetic. It is well known that a point is minimal if and only if it is almost periodic (\cite[Theorem 1.7]{Aus}).

\medskip
Let $k\geq 2$ and  $(X_1,G), \ldots ,(X_k, G)$ be $G$-systems.  Then there is a natural diagonal action of $G$ on the product space $X_1\times\cdots \times X_k$ by $g(x_1,\ldots, x_k)=(gx_1,\ldots, gx_k)$ for each $g\in G$ and $(x_1,\ldots,x_k)\in X_1\times\cdots \times X_k$, and we denote this product system by $(X_1\times\cdots\times X_k, G)$.
A $G$-system $(X,G)$ is $k$-{\it transitive} for some $k\geq 1$, if the product system $(X^{k}, G)$ is transitive. We also say $2$-transitive is {\it weakly mixing} and $1$-transitive is just transitive. We say that $(X,G)$ is $\infty$-{\it transitive}
 if $(X,G)$ is $k$-transitive for any $k\geq 1$.

\medskip
Let $(X,G)$ and $(Y,G)$ be two $G$-systems. They are {\it weakly disjoint} if the product system $(X\times Y, G)$ is transitive. A {\it joining} of $(X,G)$ and $(Y,G)$  is a nonempty $G$-invariant closed subset of $X\times Y$ that has full projection onto each coordinate. Clearly, $X\times Y$ is a joining of $(X,G)$ and $(Y,G)$. If $X\times Y$ is the only joining of $(X,G)$ and $(Y,G)$, then we say that they are {\it disjoint} and denote by of $(X,G)\perp(Y,G)$.

\medskip
Let $H$ be a subgroup of $G$, then there is also an action of $H$ on $X$ and we denote this system by $(X,H)$. A $G$-system $(X,G)$ is {\it totally transitive} if for any subgroup $H$ of $G$ with finite index, the system $(X,H)$ is transitive.

\medskip
We say a system $(Y, G)$ is a factor of a system $(X,G)$ if there is  continuous onto map $\pi: X\rightarrow Y$ such that for each $x\in X$ and $g\in G$, we always have $g\pi(x)=\pi(gx)$; in this case, we also say $(X, G)$ is an extension of $(Y,G)$.

\medskip
Let $(X,G)$ be a $G$-system. Two points $x, y\in X$ are {\it proximal} if $\inf_{g\in G}\rho(gx,gy)=0$ , where $\rho$ is a metric on $X$, and in this case we say $(x,y)$ is a {\it proximal pair}. The system $(X, G)$ is {\it proximal} if any two distinct points of $X$ are proximal and is {\it distal} if any two distinct points of $X$ are not proximal. A point $x\in X$ is {\it distal} is for any $y\in \overline{Gx}\setminus\{x\}$ , $(x,y)$ is not a proximal pair. For $x\in X$, the {\it proximal cell} of $x$ is the set of points that are proximal with $x$ and denoted by ${\rm Prox}_{G}(x)$ or ${\rm Prox}(x)$ for short.

\medskip
Given a group $G$, the Bernoulli shift over $G$  is defined to be $2^{G}=\{0,1\}^{G}$ endowed with the product topology and the action of $G$ on $2^{G}$ is defined by right translations, i.e., for $x\in 2^{G}$ and $g\in G$, $g\cdot x\in 2^{G}$ is defined by $(g\cdot x)(h)=x(hg)$ for any $h\in G$. For a subset $S\subset G$, let ${\bf 1}_{S}\in 2^{G}$ be the characteristic function of $S$, i.e.
\[{\bf 1}_{S}(g)=\begin{cases} 1& \text{ if } g\in S\\ 0 &\text{ otherwise}.\end{cases}\]

\medskip
Finally, the following notion will be frequently used. For a $G$-system $(X,G)$ and two nonempty  subsets $A,B\subset X$, define
\[ N_{G}(A, B)=\{g\in G: gA\cap B\neq \emptyset\}\]
 and when $A=\{x\}$, denote $N_{G}(\{x\},B)$ by $N_{G}(x, B)$ for short. If there is no confusion in the context, we also denote $N_{G}(A,B)$ by $N(A, B)$ for short.

\subsection{Some basic lemmas for disjointness}

The following lemma is shown in \cite[Proposition 1.1]{HY} for $\mathbb{Z}$-systems. However, the proof for general cases is similar and we reprove it here for the reader's convenience.
 \begin{lem}\label{factor}
Let $(X,G), (Y,G)$ be two $G$-systems and $(Z,G)$ be a factor of $(X,G)$. If $(X,G)$ is disjoint with $(Y,G)$ then $(Z,G)$ is also disjoint with $(Y,G)$.
\end{lem}
\begin{proof}
Let $\pi: (X,G)\rightarrow (Z, G)$ be the factor map. Let $J$ be a joining of $(Z, G)$ and $(Y,G)$. Then $\widetilde{J}=\{(x,y)\in X\times Y: (\pi(x), y)\in J\}$ is a joining of $(X, G)$ and $(Y,G)$ and hence $\widetilde{J}=X\times Y$ by the disjointness. Now the disjointness of $(Z, G)$ and $(Y,G)$ follows from that $J=\pi\times id(\widetilde{J})=Z\times Y$.
\end{proof}

A system $(X,G)$ is called an {\it $M$-system} if it is topologically transitive and the set of minimal points is dense in $X$. The following lemma is proved in \cite[Lemma 2.1]{HY} for a single map and in \cite[Theorem 3.1]{WCF} for general semigroup actions.

\begin{lem}\label{M-sys}
Let $(X,G)$ be a transitive system. Then the following assertions are equivalent.
\begin{itemize}
\item[(1)] $(X,G)$ is an $M$-system;
\item[(2)] $N(x, U)$ is a right piecewise syndetic set for any transitive point $x\in X$ and any open neighborhood $U$ of $x$.
\end{itemize}
\end{lem}


\begin{lem}\cite[Theorem 1]{CKN}
Every weakly mixing $G$-system is totally transitive. If $G$ is abelian, then every weakly mixing $G$-system is $\infty$-transitive.
\end{lem}

The following lemma is shown in \cite[Theorem 2.9]{AG} for $\mathbb{Z}$-actions but the proof for general group actions is the same. For the reader's convenience, we afford a proof here.
\begin{lem}\label{weakly disjoint}
Let $G$ be a discrete group and $(X,G)$ be a transitive system. If $(X,G)$ is weakly disjoint with all minimal $G$-systems, then $(X,G)$ is weakly disjoint with all transitive $G$-systems with dense minimal points.
\end{lem}
\begin{proof}
Let $(Y,G)$ be an arbitrary  transitive $G$-system with dense minimal points. We need to show the diagonal action of $G$ on $X\times Y$ is transitive. By the definition of transitivity, we need to show that for any nonempty open sets $A, B\subset X$ and $C,D\subset Y$, we have $N(A\times C,  B\times D)\neq \emptyset$.

Fix  nonempty open sets $A, B\subset X$ and $C,D\subset Y$. Since $(Y, G)$ is transitive, there is some $g_0\in G$ such that $C\cap g_0^{-1}D\neq\emptyset$. Thus $C\cap g_0^{-1}D$ is a nonempty open set in $Y$. Take a minimal point $y\in C\cap g_0^{-1}D$ and let $Z=\overline{G\cdot y}$. Then we have $C\cap Z\neq\emptyset$ and $D\cap Z\neq\emptyset$. It follows that $C\cap Z$ and $D\cap Z$ are nonempty open subsets of $Z$ under the subspace topology. Recall that $(X,G)$ is weakly disjoint with $(Z,G)$, since $(Z,G)$ is a minimal system. Thus $N(A\times (C\cap Z), B\times (D\cap Z))\neq \emptyset$. In particular, we have $N(A\times B, C\times D)\neq \emptyset$, and hence we show that $(X,G)$ is weakly disjoint with $(Y,G)$.
\end{proof}

A $G$-system $(X,G)$ is {\it incontractible} if, for every $n\geq 1$, the minimal points in $X^{n}$ are dense. The following lemma is shown in \cite[Theorem 4.2]{GLA} for minimal systems and the necessity is shown in \cite[Proposition 4.4]{GTWZ} for general systems.  For the completeness, we will show the other direction in the appendix.
\begin{lem}\label{incontractible}
Let $(X,G)$ be a $G$-system with dense minimal points. Then $(X,G)$ is incontractible if and only if $(X,G)$ is disjoint with all minimal proximal $G$-systems.
\end{lem}

The following result is shown in \cite[Theorem 3.8]{AK} for $\mathbb{Z}$-actions under the weakly mixing condition, which is equivalent to $\infty$-transitivity for abelian groups. Recall that a subset in a topological space is a {\it residual set} if it contains a dense $G_{\delta}$ set.
\begin{lem}\label{proximal cell}
Let $(X,G)$ be an $\infty$-transitive system. Then for every $x\in X$, the proximal cell ${\rm Prox}_{G}(x)$  is residual in $X$.
\end{lem}

\begin{proof}
For each $\varepsilon>0$, let $\Delta_{\varepsilon}=\{(x,y)\in X\times X: d(x,y)<\varepsilon\}$. Then the proximal relation of $(X,G)$ is ${\rm Prox}_{G}=\bigcap_{n=1}^{\infty}(G\Delta_{1/n})$. Thus ${\rm Prox}_{G}$ is a $G_{\delta}$ subset of $X\times X$. Hence, for each $x\in X$, ${\rm Prox}_{G}(x)$  is a $G_{\delta}$ subset of $X$. So it remains to show ${\rm Prox}_{G}(x)$  is dense in $X$ for every $x\in X$.

\medskip
Fix an $x\in X$ and a nonempty open set $U$ in $X$. It suffices to show that there is some point $y\in U$ that is proximal to $x$. Let
\[ \mathcal{F}=\{N_{G}(V, W):~~ V, W \text{ are nonempty open sets in } X\}.\]
It follows from the $\infty$-transitivity of $(X,G)$ that the family $\mathcal{F}$ has the  finite intersection property, i.e., for any $k\geq 1$ and any nonempty open sets $V_1,\ldots, V_k, W_1,\ldots, W_k$ in $X$,
\[ N_{G}(V_1, W_1)\cap \cdots\cap N_{G}(V_k, W_k)=N_{G}(V_1\times\cdots\times V_{k}, W_1\times\cdots\times W_{k})\neq \emptyset.\]
Thus, by the compactness, we have $\bigcap_{F\in \mathcal{F}}\overline{F\cdot x}\neq \emptyset$. Take some $z\in \bigcap_{F\in \mathcal{F}}\overline{F\cdot x}$ and we have $ N_{G}(x, W)\cap F\neq\emptyset$, for every open neighborhood $W$ of $z$ and every $F\in \mathcal{F}$. Let $W_{k}=B(z, 1/k)$ for each $k\geq 1$.  Now set $U_0=U$ and inductively define nonempty open sets $U_{k}$ and $g_{k}\in G$ as follows. By noting that $N_{G}(U_{k-1}, W_{k})\cap  N_{G}(x, W_{k})\neq\emptyset$, we can choose a nonempty open set $U_k$ with $\overline{U_{k}}\subset U_{k-1}$ and $g_k\in G$ such that $g_k(x)\in W_{k}$ and $g_{k}(\overline{U_{k}})\subset W_{k}$.  It follows that $\bigcap_{k=1}^{\infty}\overline{U_{k}}=\bigcap_{k=1}^{\infty} U_{k}\neq\emptyset$, and for any $y\in \bigcap_{k=1}^{\infty} U_{k}$, we have $\lim_{k\rightarrow\infty}d(g_{k}x, g_{k}y)=0$. Thus $y\in U$ and $y$ is proximal to $x$.
\end{proof}

\section{Transitive system disjoint with all minimal systems}\label{char of dis}

In this section, we give some criteria for disjointness via the notion of $m$-set  introduced by Huang and Ye in \cite{HY} and the notion of dense orbit set introduced by Glasner, Tsankov, Weiss and Zucker in \cite{GTWZ}.  Further, the relation between these two notions is also discussed.

\subsection{Characterization of disjointness via $m$-set}
Let $G$ be a discrete group and let $\mathcal{M}(G)^{\perp}$ denote the collection of  $G$-systems that are disjoint with all minimal $G$-systems.

 \medskip
The following notion is introduced in \cite{HY}.
\begin{defn}
A subset $A$ of a discrete group $G$ is an $m$-set of $G$ if there is a minimal $G$-system $(X,G)$, $x\in X$,  and a nonempty open set $U\subset X$ such that $$A\supset N_{G}(x, U)=\{g\in G: gx\in U\}.$$
\end{defn}

Let $\mathcal{P}$ be a collection of subsets of the group $G$.  A subset $\mathcal{F}$ of $\mathcal{P}$ is a {\it family} if it is hereditarily upward, i.e., for any $F_1\in \mathcal{F}$, $F_1\subset F_2$ implies $F_2\in\mathcal{F}$. For a subset $\mathcal{A}$ of $\mathcal{P}$, the family generated by $\mathcal{A}$ is defined as
$ \{ B\in \mathcal{P}: \exists A\in \mathcal{A} \text{ with } A\subset B\}$.

Clearly, all $m$-sets of $G$ form a family.  The following characterization of the family of $m$-sets of $\mathbb{Z}$ is showed in \cite[Proposition2.3]{HY}. However, it is similar for general groups.
\begin{prop}\label{m-set by min pt}
The family of $m$-sets of a countable group $G$ is generated by the subsets of $G$ whose indicator function is a minimal point in the Bernoulli shift $2^{G}$.
\end{prop}

Let $\mathcal{F}$ be a family of $\mathcal{P}$. The {\it dual family} $\mathcal{F}^{*}$of $\mathcal{F}$ is defined as
\[ \mathcal{F}^{*} =\{E\in \mathcal{P}:~~\forall F\in\mathcal{F},  E\cap F\neq \emptyset\}.\]

The following characterization of disjointness via $m$-set is also showed in \cite{HY} for $\mathbb{Z}$-actions.
\begin{lem}\label{char by m-set}
Let $(X,G)$ be a transitive $G$-system and $x\in Tran_{G}(X)$. Then
\begin{itemize}
\item[(1)]  $(X,G)\in\mathcal{M}(G)^{\perp}$ if and only if $N_{G}(x, U)\cap A\neq\emptyset$ for any open neighborhood $U$ of $x$ and  $m$-set $A$ of $G$;
\item[(2)]  $(X,G)$ is disjoint with a minimal $G$-system $(Y,G)$ if and only if $N_{G}(x, U)\cap N_{G}(y, V)\neq\emptyset$ for any open neighborhood $U$ of $x$, any point $y\in Y$ and any nonempty open set $V\subset Y$.
\end{itemize}
\end{lem}
\begin{proof}
We just show (1), as (2) can be easily followed from the proof of (1).

\medskip
Suppose $(X,G)\in \mathcal{M}(G)^{\perp}$.  Let $A$ be an $m$-set and $(Y,G)$ be an associated minimal $G$-system so that $A\supset N_{G}(y, V)$ for some $y\in Y$ and a nonempty open set $V\subset Y$. Let $J=\overline{G\cdot(x,y)}$. Then $J$ is joining since $x$ is transitive and $y$ is minimal. By the disjointness, we have $J=X\times Y$ and hence for any open neighborhood $U$ of $x$, we have
\[ N_{G}(x, U)\cap A\supset  N_{G}(x, U)\cap N_{G}(y, V)=N_{G}((x,y), U\times V)\neq \emptyset.\]

To show the contrary, let $(Y,G)$ be a minimal $G$-system and  $J$ be a joining of $X$ and $Y$. Since $J$ projects fully to each coordinate, there is $y\in Y$ with $(x,y)\in J$. Then for any  open neighborhood $U$ of $x$ and any open set $V\subset Y$, we have $N_{G}(x, U)\cap N_{G}(y, V)=N_{G}((x,y), U\times V)\neq \emptyset$, since $N_G(y, V)$ is an $m$-set. Thus $\overline{G\cdot(x,y)}\cap (U\times V)\neq\emptyset$. By noting that $\overline{G\cdot(x,y)}$ is $G$-invariant and the orbits of $x$ and $y$ under $G$ are dense in $X$ and $Y$ respectively,   we have $\overline{G\cdot(x,y)}=X\times Y$ and hence $J=X\times Y$. Thus $(X,G)$ is disjoint with $(Y,G)$, whence $(X,G)\in \mathcal{M}(G)^{\perp}$.
\end{proof}

\subsection{Characterization of disjointness via dense  orbit set}
The following notion of dense orbit set is introduced in \cite{GTWZ} and the characterization of disjointness is also implicitly given there.
\begin{defn}
Given a group $G$, a subset $A\subset G$ is a {\it dense orbit set} of $G$,  if for any minimal $G$-system $(X,G)$ and any $x\in X$, the set $A\cdot x=\{gx: g\in A\}$ is dense in $X$. Let $DO(G)$ denote the set of dense orbit sets of $G$.
\end{defn}

\begin{lem}
Let $(X,G)$ be a transitive $G$-system. Then $(X,G)\in\mathcal{M}(G)^{\perp}$   if and only if there is a point  $x\in Tran_{G}(X)$ such that for any nonempty open set $U\subset X$, $N(x, U)\in DO(G)$.

When the latter condition holds for some transitive point then it holds for every transitive point.
\end{lem}

\begin{proof}
Suppose $(X,G)\in\mathcal{M}(G)^{\perp}$. Let $(Y,G)$ be a minimal $G$-system and $y\in Y$. Let $U\subset X$ be a nonempty open set  and $x\in Tran_{G}(x)$. Let $J=\overline{G\cdot(x, y)}$. Clearly, $J$ is a joining of $X$ and $Y$. Thus $J=X\times Y$ since $(X,G)\perp (Y,G)$. It follows that $N(x, U)\cdot (x,y)\subset U\times Y$ is dense and hence $N(x,U)\cdot y$ is dense in $Y$. So $N(x, U)\in DO(G)$.

\medskip
Now suppose $(X,G)$ satisfies that for any nonempty open set $U\subset X$ and any $x\in Tran_{G}(X)$, $N(x, U)\in DO(G)$.  Fix a nonempty open set $U\subset X$. Let $(Y,G)$ be a minimal $G$-system and $y\in Y$. Let $J\subset X\times Y$ be a joining of $X$ and $Y$. Let $U\subset X$ and $V\subset Y$ be any nonempty open sets. Take a point $x_0\in Tran_{G}(X)$. Let $y_0\in Y$ be such that $(x_0,y_0)\in J$. By assumption, $N(x_0, U)\cdot y_0$ is dense in $Y$. Thus there is $g\in N(x_0,U)\subset G$ such that $g\cdot (x_0, y_0)\in U\times V$. Since $J$ is  invariant under $G$, we have $J\cap (U\times V)\neq\emptyset$. Thus  $J=X\times Y$, since $J$ is closed and $U,V$ are arbitrary. So $(X,G)\in\mathcal{M}(G)^{\perp}$.
\end{proof}

\subsection{Relations among different notions}

 Actually, the systems considered in \cite{GTWZ} have no requirement on the metrizability. It was asked in \cite[Question 9.6]{GTWZ} how to characterize dense orbit sets in countable discrete groups. This was answered in \cite{KRS} by the notion of symmetrically syndetic sets.
\begin{defn}\label{defn of ss}
Given a group $G$, a subset $A\subset G$ is  {\it symmetrically right syndetic}, if for any finite subsets $F_1\subset A$ and $F_2\subset A^{c}$,  the set $(\cap_{f_1\in F_1}f_1^{-1}A)\cap (\cap_{f_2\in F_2}f_2^{-1}A^{c})$ is right syndetic. Let $SS(G)$ denote the set of symmetrically right syndetic subset of $G$.
\end{defn}

A subset $A\subset G$ is a {\it strong dense orbit set} of $G$,  if for any minimal $G$-system $(X,G)$ with $X$ a compact Hausdorff space and any $x\in X$, the set $A\cdot x=\{gx: g\in A\}$ is dense in $X$. Let $SDO(G)$ denote the set of strong dense orbit set of $G$. The following characterization is given in  \cite{KRS}.

\begin{prop}\label{char of SDO}
Let $G$ be a discrete group. Then $$SDO(G)=SS(G)^{*}=\{A\subset G: ~\forall B\in SS(G), ~~A\cap B\neq \emptyset\}.$$
\end{prop}

Similarly, we can also characterize the dense orbit set via the notion of $m$-set. For a group $G$, we use $m(G)$ to denote the family of $m$-sets of $G$.
\begin{prop}\label{char of DO}
Let $G$ be a countable discrete group. Then $$DO(G)=m(G)^{*}=\{A\subset G: ~\forall B\in m(G), ~~A\cap B\neq \emptyset\}.$$
\end{prop}
\begin{proof}
Suppose that $A\subset G$ is a dense orbit set and $B\subset G$ is an $m$-set. Then it follows from the definition of $m$-set that there is a minimal $G$-system $(X, G)$ with a point $x\in X$ and a nonempty open set $U\subset X$ such that $B\supset N(x, U)$. While we have $Ax\cap U\neq\emptyset$, since $A$ is a  dense orbit set. Thus $A\cap B\neq \emptyset$ and hence $DO(G)\subset m(G)^{*}$.

\medskip
Assume $A$ is in $m(G)^{*}$. Then each $m$-set has nonempty intersection with $A$. In particular, for any minimal $G$-system $(X, G)$, any $x\in X$ and any nonempty open set $U\subset X$, we have $A\cap N(x, U)\neq\emptyset$, which implies that $Ax$ is dense in $X$. Thus $A$ is a dense orbit set and it follows that $DO(G)\supset m(G)^{*}$.

\medskip
Therefore, we have $DO(G)=m(G)^{*}$.
\end{proof}

The following equivalent characterization of symmetrically right syndetic set in given in \cite{KRS} via the set of recurrence time in the universal minimal flow.
\begin{prop}\cite[Proposition 6.11]{KRS}\label{char of ss}
Let $G$ be a discrete group. A subset $A\subset G$ is symmetrically syndetic if and only if the universal minimal flow contains an clopen subset $U$ and a point $x$ such $A=N(x, U)$.
\end{prop}

We can also give a characterization of symmetrically right syndetic set via the minimal points in the Bernoulli shift.
\begin{prop}\label{char of ss via min pt}
Let $G$ be a discrete group. Then a nonempty proper subset $A\subset G$ is symmetrically right syndetic if and only if ${\bf 1}_{A}$ is a minimal point in $2^{G}$.
\end{prop}
\begin{proof}
Note that for any finite subset $F\subset G$ and a partial function $\varphi: F\rightarrow \{0,1\}$, we have
\begin{align*}
 N({\bf 1}_{A}, U[\varphi])&=\{g\in G: g{\bf 1}_{A}\mid_{F}=\varphi\}\\
&=\{g\in G: {\bf 1}_{A}(fg)=\varphi(f), \forall f\in F\}\\
&=\{g\in G: f_1g\in A, f_2g\in A^{c}, \forall f_1\in F_1 , f_2\in F_2 \}\\
&=(\cap_{f_1\in F_1}f_1^{-1}A )\cap(\cap_{f_2\in F_2}f_2^{-1}A^{c}),
\end{align*}
where $F_1=\varphi^{-1}(1), F_2=\varphi^{-1}(0)$ and $U[\varphi]=\{ \xi\in 2^{G}: \xi\mid_{F}=\varphi\}$.

\medskip
Now assume that $A\subset G$ is a symmetrically right syndetic subset of $G$. Then $ N({\bf 1}_{A}, U[\varphi])$ is right syndetic for each finite subset $F\subset G$ and a partial function $\varphi: F\rightarrow \{0,1\}$. Thus ${\bf 1}_{A}$ is a minimal point in $2^{G}$.

\medskip
Conversely, assume $A\subset G$ is a nonempty proper subset of $G$ and ${\bf 1}_{A}$ is a minimal point in $2^{G}$. Then for  each finite subset $F\subset G$ and a partial function $\varphi: F\rightarrow \{0,1\}$, the set $(\cap_{f_1\in F_1}f_1^{-1}A )\cap(\cap_{f_2\in F_2}f_2^{-1}A^{c})=N({\bf 1}_{A}, U[\varphi])$ is right syndetic and hence $A$ is symmetrically right syndetic.
\end{proof}

Now, combining Proposition \ref{m-set by min pt}, \ref{char of SDO},  \ref{char of DO} and \ref{char of ss via min pt}, we can establish the relations among the notions of $m$-set, strong dense orbit set and symmetrically right syndetic set as follows. Thus the question of the authors in \cite[Question 9.6]{GTWZ} can be answered via the notions of $m$-set, symmetrically syndetic set and the minimal points in the Bernoulli shift.

\begin{prop}
Let $G$ be a countable group. Then the family  of $m$-sets of $G$ is generated by the collection   of  symmetrically right syndetic subsets of $G$ and one has $SDO(G)=SS(G)^{*}, DO(G)=m(G)^{*}$.
\end{prop}

 Finally, combining Proposition \ref{char of ss via min pt} and the following lemma, we conclude that for a countable group $G$, the set of minimal points in $2^{G}$ is dense, which has been proven in \cite[Theorem 5.3.6]{GJS} and \cite[Lemma 4.2]{GTWZ}.

\begin{lem}\label{ss extends partial function}
Let $G$ be a discrete group and $F$ be a finite subset of $G$. Then for any map $\varphi: F\rightarrow \{0,1\}$, there is a symmetrically syndetic subset $A\subset G$ such that the characteristic function ${\bf 1}_{A}$ extends $\varphi$, that is ${\bf 1}_{A}\mid_{F}=\varphi$. In particular, the set $\{ {\bf 1}_A: A\subset G \text{ is symmetrically syndetic}\}$ is dense in $2^{G}$.
\end{lem}
\begin{proof}
Let $M$ be the universal minimal flow of $G$ and let $F_1=\varphi^{-1}(1)$ and $F_2=\varphi^{-1}(0)$. Take a point $x\in M$. It is well known that $M$ is a compact Hausdorff totally disconnected space (\cite{BF}).  In addition, the action of $G$ on $M$ is free by Ellis (\cite{Ellis}) or Veech (\cite{Veech}). In particular, $Fx$ consists of $|F|$ distinct points. Thus there is a clopen subset $U\subset M$ such that $F_1x\subset U$ and $F_2x\cap U=\emptyset$.  Then it follows from  Proposition \ref{char of ss} that $N(x, U)$ is a symmetrically syndetic subset of $G$. Clearly, $F_1\subset N(x, U)$ and $F_2\cap N(x,U)=\emptyset$. Let $A=N(x, U)$. Then $A$ is what we are looking for.
\end{proof}

\section{Proof of Theorem \ref{thmA}}

The following result is crucial for us but the proof of it is very technical and is given in Section 6.

\begin{prop}\label{m-set}
Let $G$ be a countably infinite group. Then every right thickly syndetic subset of $G$ contains an $m$-set.
\end{prop}

Using the above proposition, we are able to show

\begin{thm}
Let $G$ be a countable group and $(X,G)$ be a transitive system. If $(X,G)$ is disjoint with all minimal $G$-systems, then $(X,G)$ is an $M$-system without nontrivial minimal factor.
\end{thm}
\begin{proof}
If $G$ is a finite group then $(X,G)$ is also a minimal system. Since a nontrivial minimal system is never disjoint with itself,  we conclude $(X,G)$ is a trivial system when $G$ is finite and the theorem holds trivially. So we may assume that $G$ is countably infinite.

\medskip
Let $x\in X$ be a transitive point and $U$ be a nonempty open neighborhood of $x$. By Lemma \ref{char by m-set}, $N(x, U)\cap A\neq \emptyset$ for any $m$-set $A\subset G$. It follows from Theorem \ref{m-set} that $N(x, U)$ has nonempty intersection with every right thickly syndetic subset of $G$.  Thus, by Lemma \ref{piecewise synd}, $N(x, U)$  is right piecewise syndetic. By Lemma \ref{M-sys}, $(X,G)$ is an $M$-system.

\medskip
Let $(Z,G)$ be a minimal factor of $(X,G)$. By Lemma \ref{factor}, $(Z,G)$ is also disjoint with all minimal $G$-systems. In particular, $(Z,G)$ is  disjoint with itself. Thus $(Z,G)$ is trivial.
\end{proof}

\begin{defn}
Let $(X,G)$ and $(Y,G)$ be two $G$-systems, and let $D\subset X$, $E\subset Y$ be two subsets. We say $(D,E)$ is a {\it resonating pair} of $X$ and $Y$ if for any $y\in E$ and any open neighborhood $V$ of $y$, and any nonempty open set $U\subset X$ , there is some $x\in U\cap D$ such that   $N_{G}((x,y), U\times V)=\{g\in G: gx\in U, gy\in V\}$ is right syndetic.
\end{defn}

We remark that it follows from the definition that $D$ is necessary dense in $X$.

Now Theorem \ref{thmA} follows from the next theorem.
\begin{thm}\label{stronger than thmA}
Let $(X,G)$ be a transitive $G$-system. Then the following statements are equivalent.
\begin{itemize}
\item[(1)] $(X,G)$ is disjoint with all minimal $G$-systems.
\item[(2)]  $(X, G)$ is   $\infty$-transitive and for every minimal $G$-system $(Y,G)$, there is a countable set $D_{Y}\subset X$  such that   $(D_Y, Y)$ is  a resonating pair of $X$ and $Y$.
\item[(3)] $(X, G)$ is   $\infty$-transitive and there is a countable set $D\subset X$ consisting of minimal points such that for every minimal $G$-system $(Y,G)$,   $(D, Y)$ is a resonating pair of $X$ and $Y$.
\end{itemize}
\end{thm}

\begin{proof}
Clearly,   (3) $\Longrightarrow$ (2).

\medskip
(2) $\Longrightarrow$ (1)  Fix a minimal $G$-system $(Y, G)$ and let $(D_Y, E)$ be the associated resonating pair as in (2).  Let $D=D_{Y}$.

\medskip
Fix a joining $J$ for $X$ and $Y$. It follows from Lemma \ref{proximal cell} that
${\rm Prox}(D):=\bigcap_{x\in D}{\rm Prox}(x)$ is a residual set in $X$. Take $x_0\in D$. Then here is some    $y\in Y$ with $(x_0,y)\in J$. Let $U$ be any nonempty open set in $X$ and $V$ be an open neighborhood of $y$. Take a nonempty open set $W\subset X$ and $\varepsilon >0$ such that $B_{\varepsilon}(\overline{W})\subset U$. According to the assumption in (2), there is some $x\in D\cap W$ such that $N_{G}((x,y), W\times V)$ is right syndetic. Since $(x,x_0)$ is proximal, we have that $\{g\in G: d(gx,gx_0)<\varepsilon\}$ is right thick. Thus there is some $g\in G$ such that
\[ g(y)\in V,\ \ g(x)\in W\ \ \text{ and }\ \  d(gx, gx_0)<\varepsilon.\]
So $g(x_0)\in U$ and hence $(gx_0, gy)\in J\cap (U\times V)$. Since $U$ and $V$ are arbitrary, we have $X\times \{y\}\subset J$. It now follows from the minimality of $(Y,G)$ that $J=X\times Y$, and hence $(X,G)$ is disjoint with $(Y,G)$.

\medskip
Since $(Y,G)$ is arbitrary, we conclude that $(X,G)$ is disjoint with all minimal $G$-systems.

\medskip
(1) $\Longrightarrow$ (3) Suppose that $(X, G)$ is disjoint with all minimal $G$-systems. It follows from Proposition \ref{m-set} that the set of minimal points in $(X,G)$  is dense.  By Lemma \ref{weakly disjoint}, $(X,G)$ is weakly disjoint with itself and hence $(X,G)$ is $2$-transitive. Assume that  $(X,G)$ is $k$-transitive for some $k\geq 2$. Then Lemma  \ref{incontractible} implies that $(X^{k},G)$ has dense minimal points. Thus, by Lemma \ref{weakly disjoint},   $(X,G)$ is weakly disjoint with $(X^{k}, G)$ and hence $(X,G)$ is $(k+1)$-transitive. Therefore, $(X,G)$ is $\infty$-transitive.

\medskip
Let $D$ be any countable dense subset of $X$ consisting of minimal points. We are going to show that such $D$ satisfies the requirements. To the contrary, assume that the condition in (3) does not hold for $D$. Then there is a minimal $G$-system $(Y,G)$ with a point $y\in Y$ and its open neighborhood $V$, and there is a nonempty open set $U\subset X$ such that $N_{G}((x,y), U\times V)$ is not right syndetic for any point $x\in D\cap U$.

\medskip
\noindent {\bf Claim}. For every $x\in D\cap U$, there is a $G$-invariant subset $J_{(x,y)}\subset (U\times V)^{c}$ such that $x\in {\rm Proj}_{1}(J_{(x,y)})$, where ${\rm Proj}_{1}$ is the projection of $X\times Y$ to the first coordinate.

\medskip
\begin{proof}[Proof of Claim]
Since $N_{G}((x,y), U\times V)$ is not right syndetic, the set $A=N_{G}((x,y), (U\times V)^{c})$ is right thick.   As $x$ is a minimal point, the set $N_{G}(x, B_{\varepsilon}(x))$ is right syndetic for every $\varepsilon>0$. Thus, for any $\varepsilon>0$ and finite subset $F\subset G$, there is some $h\in N_{G}(x, B_{\varepsilon}(x))$ with $Fh\subset A$,  since  the set $\{g\in G: Fg\subset A\}$ is right thick by Lemma \ref{thick translate}. Let $(F_n)_{n=1}^{\infty}$ be an increasing sequence of finite subsets of $G$ such that $G=\bigcup_{n=1}^{\infty}F_n$. Then we can choose a sequence $(g_k)$ in $G$ such that
\[ g_{k}\in N_{G}(x, B_{1/k}(x)) \text{ for each } k, \ \ \text{ and }\ \ F_{k}g_{k}\subset A.\]
By passing to some subsequence, we may assume that
\[ \lim_{k\rightarrow \infty} g_{k}\cdot (x,y)=(x_0, y_0).\]
It is clear that $x_0=x$, since $g_{k}\in N_{G}(x, B_{1/k}(x))$. For any $g\in G$, there is some  $k_0\geq 1$ such that $g\in F_{k_0}$ and then we have
\[ g\cdot (x_0,y_0)=g( \lim_{k\rightarrow \infty} g_{k}\cdot (x,y))= \lim_{k\rightarrow \infty} gg_{k}\cdot (x,y)\in (U\times V)^{c},\]
since $F_{k_0}g_{k}\subset F_{k}g_{k}\subset A$ for any sufficiently large $k$. Hence $\overline{G\cdot (x_0, y_0)}\subset (U\times V)^{c} $. Now set $J_{(x,y)}=\overline{G\cdot (x_0, y_0)}$. Then $J_{(x,y)}\subset (U\times V)^{c} $ and $x=x_0\in{\rm Proj}_{1}(J_{(x,y)})$. This shows the claim.
\end{proof}

\medskip
Now set \[ J=\overline{\bigcup_{x\in D\cap U}J_{(x,y)}}.\]
It follows from the above claim that $J\subset (U\times V)^{c},\  D\cap U\subset {\rm Proj}_{1}(J)$ and $J$ is invariant under the diagonal action of $G$. Since $D$ is dense in $X$ and ${\rm Proj}_{1}(J)$ is closed, we have $U\subset {\rm Proj}_{1}(J)$. Since $(X,G)$ is transitive and ${\rm Proj}_{1}(J)$ is a $G$-invariant closed subset, we have ${\rm Proj}_{1}(J)=X$. In addition,  it follows from the minimality of $(Y,G)$ that  $J$ has full projection onto $Y$. Thus $J$ is a joining of $(X,G)$ and $(Y,G)$. Hence the disjointness implies that $J=X\times Y$. But this contradicts that $J\subset (U\times V)^c$.
\end{proof}

\medskip
 To show that the condition in Theorem \ref{stronger than thmA}  can be verified for a concrete system, next we construct the  resonating pair for  the Bernoulli shift $2^{G}$ of a countable group $G$. However, our construction depends on the separated covering property for minimal systems introduced in \cite{GTWZ}. Note that we are not aiming to reprove their result but  to understand it from our point of view.

Recall that a minimal $G$-system $(X,G)$ has {\it separated covering property} ({\bf SCP} for short) if for any nonempty open set $U\subset X$ and any finite subset $F\subset G$, there are $g_1,\ldots, g_n\in G$ for some $n\in\mathbb{N}$ such that
\begin{enumerate}
\item $Fg_i\cap Fg_{j} =\emptyset, \forall i\neq j$ and
\item $X=g_1U\cup \cdots\cup g_{n}U$.
\end{enumerate}

The following lemma is a direct consequence of {\bf SCP}.
\begin{lem}\label{sep+syn}
Let $G$ be a countably infinite group and $(X,G)$ be a minimal $G$-system. If $(X,G)$ has {\bf SCP}, then for any nonempty open set $U\subset X$ and any finite subset $F\subset G$, there are finitely many $g_1,\ldots, g_n\in G$ such that
\begin{enumerate}
\item $Fg_{i}\cap Fg_{j}=\emptyset, \forall i\neq j\in\{1,\ldots,n\}$ and
\item $G=g_1N(x, U)\cup g_2N(x, U)\cup\cdots\cup g_{n}N(x,U)$ for any $x\in U$.
\end{enumerate}
\end{lem}
\begin{proof}
Suppose that $g_1,\ldots,g_n\in G$ satisfy that $Fg_i\cap Fg_{j} =\emptyset, \forall i\neq j$ and  $X=g_1U\cup \cdots\cup g_{n}U$. Then for any $g\in G$ and $x\in U$, there is some $g_i$ such that $gx\in g_{i}U$ and hence $g_{i}^{-1}g\in N(x, U)$. Thus $g\in  g_i N(x, U)$. This shows that $G=g_1N(x, U)\cup g_2N(x, U)\cup\cdots\cup g_{n}N(x,U)$ for any $x\in U$.
\end{proof}

In \cite{GTWZ}, the authors show that a minimal $G$-system is disjoint from the Bernoulli shift if and only if it has {\bf SCP} and they establish that every minimal $G$-system has {\bf SCP} for any infinite discrete group $G$. Next, we verify the condition in Theorem \ref{stronger than thmA} for the Bernoulli shift under the assumption that every minimal $G$-system has {\bf SCP}. We restate this as the following corollary. 

\begin{cor}\label{Bernoulli disjoint}
For a countable infinite group $G$, the Bernoulli shift $2^{G}$ is disjoint from all minimal $G$-system.
\end{cor}
\begin{proof}

First note that $2^{G}$ is $\infty$-transitive (see \cite[Lemma 2.3]{GTWZ} where it was shown that $2^{G}$ is transitive and similarly, $(2^{G})^{k}\cong (2^{k})^{G}$ is transitive for any $k\geq 2$). According to Theorem \ref{stronger than thmA}, it suffices to show that for each minimal $G$-system $(Y, G)$ there is a resonating pair $(D_Y, Y)$ of $2^{G}$ and $Y$, where $D_{Y}\subset 2^{G}$ is a countable dense subset of $2^{G}$.

\medskip
Fix a minimal $G$-system $(Y, G)$.

\medskip
\noindent{\bf Claim}. For any nonempty open sets $U\subset 2^{G}$ and $V\subset Y$ in $X$ and $Y$ respectively, there is $x_{U,V}\in U$ such that for any $y\in V$, $N(x_{U,V}, U)\cap N(y, V)$ is right syndetic.
\begin{proof}[Proof of Claim]
Take a finite subset $F\subset G$ and $\varphi: F\rightarrow \{0,1\}$ such that $U[\varphi]=\{x\in 2^{G}: x|_{F}=\varphi\}\subset U$. By Lemma \ref{sep+syn}, there are $f_1,\ldots, f_{n}\in G$ such that
\begin{enumerate}
\item $Ff_{i}\cap Ff_{j}=\emptyset, \forall i\neq j\in\{1,\ldots,n\}$ and
\item $G=f_1N(y, V)\cup \cdots\cup f_{n}N(y, V)$ for any $y\in V$.
\end{enumerate}
Set $E=\bigcup_{i=1}^{n}Ff_{i}^{-1}$. Now it follows from Lemma \ref{separated syndetic} that there is a right syndetic set $A\subset G$ with $e\in A$ such that for any $a\neq a'\in A$, $Ea\cap Ea'=\emptyset$. We define  $x_{U,V}\in 2^{G}$ by
\[ x_{U,V}(g)=\begin{cases}
\varphi(f), & g\in EA \text{ with } g=ff_i^{-1}a \text{ for some } f\in F, i\in\{1,\ldots, n\},\\
0, & g\in G\setminus(EA).
\end{cases}\]

\begin{figure}
\centering
 \begin{tikzpicture}
 \tikzset{
    vertex/.style={circle,draw,minimum size=3em},
    edge/.style={->,> = latex'}
}
\draw (-4.5,-7) rectangle (8.5, 2.5);

\draw (-4,0) rectangle (8, 2);
\node[vertex] (1) at (-2,1) {$~~Ff_{1}^{-1}~~$ };
\node[vertex] (1) at (1,1) {$~~Ff_{2}^{-1}~~$ };
\node  at (3.5,1)  {$\cdots$ };
\node[vertex] (1) at (6,1) {$~~Ff_{n}^{-1}~~$ };

\draw (-4,-3) rectangle (8, -1);
\node[vertex] (1) at (-2,-2) {$Ff_{1}^{-1}a_1$ };
\node[vertex] (1) at (1,-2) {$Ff_{2}^{-1}a_1$ };
\node  at (3.5,-2)  {$\cdots$ };
\node[vertex] (1) at (6,-2) {$Ff_{n}^{-1}a_1$ };

\draw (-4,-6) rectangle (8, -4);
\node[vertex] (1) at (-2,-5) {$Ff_{1}^{-1}a_2$ };
\node[vertex] (1) at (1,-5) {$Ff_{2}^{-1}a_2$ };
\node  at (3.5,-5)  {$\cdots$ };
\node[vertex] (1) at (6,-5) {$Ff_{n}^{-1}a_2$ };

\node at (2, -6.5) {$\cdots$};
   \end{tikzpicture}

\end{figure}

Since $e\in A$, we have $x_{U,V}|_{F}=\varphi$ and hence $x_{U,V}\in U[\varphi]\subset U$. 

Next we can show that $N(x_{U,V}, U)\cap N(y, V)$ is right syndetic for any $y\in V$. Set $S_{y}=N(x_{U,V}, U)\cap N(y, V)$ for each $y\in V$. Fix
$y\in V$. Since $A$ is right syndetic, there are $g_1,\ldots, g_{m}\in G$ such that $G=g_1A\cup\cdots\cup g_{m}A$. Given $g\in G$, there is some $i\in\{1,\ldots, m\}$ such that $g_{i}^{-1}g\in A$. Further, there is some $j\in \{1,\ldots,n\}$ such that $f_{j}^{-1}g_{i}^{-1}g\in N(y, V)$. On the other hand, for any $f\in F$,
\[f_{j}^{-1}g_{i}^{-1}g x_{U,V}(f)=x_{U,V}(ff_{j}^{-1}g_{i}^{-1}g)=\varphi(f).\]
Thus $f_{j}^{-1}g_{i}^{-1}g\in N(x_{U,V}, U[\varphi])\subset N(x_{U,V}, U)$. This shows that $G=\bigcup_{i=1}^{m}\bigcup_{j=1}^{n} g_{i}f_j S_{y}$. So $S_{y}$ is right syndetic.
\end{proof}

Now let $(U_i)_{i=1}^{\infty}$ and $(V_j)_{j=1}^{\infty}$ be countable basis in $2^{G}$ and $Y$ respectively. Set $D_{Y}=\{ x_{U_i, V_{j}}\in 2^{G}: i,j=1,2,\ldots\}$. Then it follows from the claim above that $(D_{Y}, Y)$ is a resonating pair for $2^{G}$ and $Y$. Thus $2^{G}$ is disjoint from $Y$ by Theorem \ref{stronger than thmA}. 
\end{proof}

We remark that the requirement of $G$ being infinite is necessary in Corollary \ref{Bernoulli disjoint}. First, for a finite group $G$, $(2^{G}, G)$ is not transitive,  let alone being $\infty$-transitive and then we can not apply Theorem \ref{stronger than thmA}.  This can be seen that the action of $G$ on the set of subsets of $G$ preserves the cardinality and hence cannot be transitive. In addition, for finite group $G$, $2^{G}$ may not be disjoint with a minimal $G$-system. For example, let $G=\mathbb{Z}/2\mathbb{Z}$.  Then $G$ acts minimally on a set $X=\{p,q\}$ by permutations. But $$J=\{(00, p), (00,q), (01,p), (10,q), (11,p),(11,q)\}$$ is a joining of $2^{G}$ and $X$, which is not equal to $2^{G}\times X$.

\medskip
Finally, we remark that both in \cite{Ber} and \cite{GTWZ} the authors show that the Bernoulli shift $2^{G}$ for any infinite group $G$ is disjoint with all minimal compact Hausdorff $G$-systems. However,  it seems at present our Theorem \ref{stronger than thmA} cannot be generalized to work in the case of uncountable groups not only due to the construction of $m$-set in a thickly syndetic set as in Proposition \ref{m-set} relies on the countability of the group but also the residual property of proximal cells in Lemma \ref{proximal cell} cannot avoid the metrizability.

\section{Proof of Theorem \ref{thmC}} 

In this section we show some sufficient conditions for a transitive $G$-system to be disjoint from all minimal $G$-systems, and Theorem \ref{thmC} follows from the results in the subsequent subsections.

\subsection{Dense small periodic sets}
\begin{defn}
Let $G$ be a discrete group. A $G$-system $(X,G)$ has {\it dense small periodic sets} if for any nonempty open subset $U$ of $X$ there is a nonempty closed subset $E\subset U$ and a finite index subgroup $H$ of $G$ such that $E$ is invariant under the action of $H$.
\end{defn}

Recall that a system $(X,G)$ is totally transitive if $(X, H)$ is transitive for any subgroup $H$ of $G$ with finite index.

\begin{lem}\label{trans for subgp}
Let $(X,G)$ be a  $G$-system having dense small periodic sets. If $(X,G)$ is totally transitive, then
\begin{itemize}
\item[(1)]  it is weakly mixing and
\item[(2)] $Tran_{G}(X)=Tran_{H}(X)$ for any finite index subgroup $H$ of $G$.
\end{itemize}
\end{lem}

\begin{proof}
(1) Let $U_1, U_2\subset X$ and $V_1,V_2\subset X$ be any nonempty open sets. We need to show there is some $g\in G$ such that $g(U_1\times V_1)\cap (U_2\times V_2)\neq \emptyset$. By the transitivity of  $(X,G)$, there is some $f\in G$ such that $W=f(U_1)\cap U_2$ is nonempty and open. Thus there is a nonempty closed subset $A\subset W$ and a finite index subgroup $H$ of $G$ such that $A$ is $H$-invariant, since $(X,G)$ has dense small periodic sets.  It follows that  $hf(U_1)\cap U_2\neq\emptyset$ for any $h\in H$, since
\[hf(U_1)\cap U_2=h(fU_1\cap h^{-1}U_2)
\supset h(A\cap h^{-1}A)=h(A)=A.\]
As $(X,G)$ is totally transitive and $f(V_1)$ is a nonempty open set, there is some $h_0\in H$ such that $h_0f(V_1)\cap V_2\neq\emptyset$. Thus, by letting $g=h_0f$, we have $g(U_1\times V_1)\cap (U_2\times V_2)\neq \emptyset$.

\medskip
(2) Let $H$ be a finite index subgroup of $X$.  By Lemma \ref{group} (3), there is a subgroup $K$ of $H$ that is normal in $G$ and has finite index in $G$. Clearly,  $Tran_{K}(X)\subset Tran_{H}(X)\subset Tran_{G}(X)$. Thus it suffices to show $Tran_{G}(X)\subset Tran_{K}(X)$.

\medskip
 To the contrary, we assume there is a point  $x\in Tran_G(X)\setminus Tran_{K}(X)$. Let $A=\overline{K\cdot x}$.    Since $\overline{G\cdot x}=X$ and $K$ has finite index in $G$, $X$ is a union of finite many homeomorphic copies of $A$ and hence $int(A)\neq\emptyset$.     Since $int(A)$ is invariant under $K$ and $K$ is normal in $G$, we have either $int(A)\cap g (int(A))=\emptyset$ or $A=gA$ for any $g\in G$. Indeed,  if $int(A)\cap g (int(A))\neq\emptyset$, then $K(int(A)\cap g (int(A)))=K(int(A))\cap g K (int(A)))=int(A)\cap g (int(A))$.  Thus $\overline{int(A)\cap g (int(A))}=A$, which implies that $A=gA$.   So  $\{g\in G: g(int(A))\cap int(A)\neq \emptyset, g(int(A))\cap (X\setminus A)\neq\emptyset\}=\emptyset$. Hence $(X,K)$ is not weakly mixing. However, it follows from Lemma \ref{group} (1) and (2) that $(X,K)$ is totally transitive and has dense small periodic sets as well. Hence $(X,K)$ is also weakly mixing.  This contradiction shows the clause.
\end{proof}

\begin{prop}\label{dsps}
Let  $(X,G)$ be a $G$-system. Then $(X,G)$ has dense small periodic sets if and only if for any nonempty open subset $U$ of $X$ there exist $p\in X$ and a finite index subgroup $H$ of $G$ such that $N_{H}(p, U)=\{h\in H: h(p)\in U\}$ is a right thick subset of $H$.
\end{prop}
\begin{proof}
Suppose that $(X,G)$ has dense small periodic sets. Then for any nonempty open subset $U$ of $X$, there is a nonempty closed subset $A\subset U$ and  a finite index subgroup $H$ of $G$ such that $A$ is $H$-invariant. Take any point $p\in A$. Then $N_{H}(p, U)=H$ is  right thick in $H$.

\medskip
Conversely, assume that for any nonempty open subset $U$ of $X$ there exist $p\in X$ and a finite index subgroup $H$ of $G$ such that $N_{H}(p, U)=\{h\in H: h(p)\in U\}$ is a right thick subset of $H$. Let $W$ be any nonempty open subset of $X$. Take a nonempty open subset $V$ with $\overline{V}\subset W$. Applying the assumption to $V$,  there is $p\in X$ and a finite index subgroup $H$ of $G$ such that $N_{H}(p, V)=\{h\in H: h(p)\in U\}$ is a right thick subset of $H$. Let $H_1\subset H_2\subset \cdots$ be an increasing sequence of finite subsets of $H$ satisfying $H=\bigcup_{n=1}^{\infty}H_n$. Now there is a sequence $(h_n)$ in $H$ such that $H_nh_n \subset N_{H}(p, V)$. Let $(n_{i})$ be an increasing subsequence such that $\lim_{i\rightarrow \infty} h_{n_i}p=z$. Let $A=\overline{H\cdot z}$.

\noindent {\bf Claim}.  $A$ is a nonempty $H$-invariant closed subset contained in $\overline{V}\subset W$.

 \medskip
 Clearly, it suffices to show that $hz\in \overline{V}$ for each $h\in H$. Fix an $h\in H$. Then for any sufficiently large $n_i$, $h\in H_{n_i}$ and $hh_{n_i}\in N_{H}(p, V)$, i.e., $hh_{n_i}p\in V$. Thus
 \[ hz=h(\lim_{i\rightarrow \infty} h_{n_i}p)= \lim_{i\rightarrow \infty} hh_{n_i}p\in \overline{V}.\]

\medskip
Therefore, $(X,G)$ has dense small periodic sets.
\end{proof}

\begin{prop}\label{char of dense small periodic set}
Let $G$ be a discrete group and $(X,G)$ be a transitive $G$-system. Then $(X,G)$ is totally transitive and has dense small periodic sets if and only if there is a transitive point $x$ satisfying condition ($\star$), that is
\begin{enumerate}
 \item[($\star$)] for any neighborhood $U$ of $x$ any $g\in G$  there is a finite index subgroup $H$ of $G$  such that $N_{H}(x, U, g)=\{h\in H: hg(x)\in U\}$ is a right thick subset of $H$.
\end{enumerate}
When the condition ($\star$) holds for some point then it will hold for every transitive point of $X$.
\end{prop}
\begin{proof}
Suppose that $(X,G)$ is totally transitive and has dense small periodic sets. Let $x\in Tran_{G}(X)$ be any transitive point. Then, by Lemma \ref{trans for subgp} (2), we have $x\in Tran_{H}(X)$ for any finite index subgroup $H$ of $G$.

By Proposition \ref{dsps}, for any neighborhood $U$ of $x$, there exist $p\in X$ and a finite index subgroup $H$ of $G$ such that $N_{H}(p, U)=\{h\in H: h(p)\in U\}$ is right thick in $H$.  Let $H_1\subset H_2\subset \cdots$ be a sequence of finite subsets of $H$ satisfying $H=\bigcup_{n=1}^{\infty}H_n$. Then there is a sequence $(h_n)$ in $H$ such that $H_nh_n\cdot p\subset U$, for each $n\geq 1$. It follows from the continuity that for each $n\geq 1$,  there is an open neighborhood $V_n$ of $p$ with $H_nh_n(V_n)\subset U$.  Now for any $g\in G$, we have $g(x)\in Tran_{H}(X)$. Thus, for each $n\geq 1$,  there exists some $h'_{n}\in H$ with $h'_{n}g(x)\in V_n$.  So we have $H_nh_nh'_ng(x)\subset U$ for each $n\geq 1$. Thus $N_{H}(x, U,g)$ is right thick in $H$.

\medskip
Conversely, assume that there is some $x\in Tran_{G}(X)$ satisfying the condition ($\star$). To show that $(X,G)$ is totally transitive and has dense small periodic sets, by Proposition \ref{dsps}, it suffices to show that for any nonempty open subset $V$ of $X$ and any finite index subgroup $H$ of $G$, $N_{H}(x, V)\neq \emptyset$ and there is some $p\in X$ and a finite index subgroup $K$ of $G$ such that $N_{K}(p, V)$ is right thick in $K$.

Now fix a nonempty open subset $V$ of $X$ and a finite index subgroup $H$ of $G$. Since $x$ is a transitive point, there is some $g\in G$ with $g(x)\in V$. Further, take a small open neighborhood $U$ of $x$ so that $g(U)\subset V$.  Applying the condition ($\star$) to $V$ and $g$, there is a finite index subgroup $L$ of $G$ such that $N_{L}(x, U, g^{-1})$ is right thick in $L$. Then for any $\ell\in  N_{L}(x, U, g^{-1})$, we have $\ell g^{-1}(x)\in U$ and hence $g\ell g^{-1}(x)\in V$, since $g(U)\subset V$.  Let $K=gLg^{-1}$ and then we have
\[ N_{K}(x, V)\supset gN_{L}(x, U, g^{-1})g^{-1},\]
which is right thick in $K$. So, by letting $p=x$, we conclude that $N_{K}(p,V)$ is right thick in $K$. By Lemma \ref{syn+thick}, we have $N_{K}(x, V)\cap H\subset N_{H}(x, V)$ is also right thick in $H$. In particular, $ N_{H}(x, V)\neq\emptyset$.
\end{proof}

\begin{lem}\label{m-set passing to subgroup}
Let $A\subset G$ be an $m$-set. Then there is some $g\in G$ such that for any finite index subgroup $H$ of $G$, the set $H\cap gA$ is an $m$-set for $H$ and hence is right syndetic in $H$.
\end{lem}
\begin{proof}
By the definition of $m$-set, there is a minimal $G$-system $(X,G)$ , $x\in X$ and a nonempty open set $U\subset X$ such that $N_{G}(x, U)\subset A$. Take an open neighborhood $V$ of $x$ and  $g\in G$ such that $g^{-1}V\subset U$. Then we have $g^{-1}N_{G}(x, V)\subset N_{G}(x, U)$ and hence $N_{G}(x, V)\subset gN_{G}(x, U)\subset gA$. Fix a finite index subgroup $H$ of $G$. Since $H$ has finite index in $G$, $x$ is also a minimal point for $H$ acting on $X$.  Thus the set $\{h\in H: hx\in V\}$ is an $m$-set and hence is right syndetic in $H$. Therefore, it follows from $\{h\in H: hx\in V\}\subset H\cap N_{G}(x, V)$ that $H\cap g A$ is an $m$-set and also right syndetic in $H$.
\end{proof}

Now we are ready to show
\begin{thm}\label{dense small periodic set}
Let $G$ be a discrete group and $(X,G)$ be a $G$-system. If $(X,G)$ is totally transitive and has dense small periodic sets, then $(X,G)\in \mathcal{M}(G)^{\perp}$.
\end{thm}
\begin{proof}
 According to Lemma \ref{char of dense small periodic set}, there is a transitive point $x$ that satisfies the property ($\star$). Thus by Lemma \ref{char by m-set}, it remains to show that for any neighborhood $U$ of $x$ and any $m$-set $A$, $N_{G}(x, U)\cap A\neq \emptyset$.

 \medskip
By Lemma \ref{m-set passing to subgroup}, there is some $g\in G$ such that $H\cap gA$ is right syndetic in $H$ for any finite index subgroup $H$ of $G$.  However,  it follows from the property ($\star$) that there is a finite index subgroup $K$ of $G$ such that $N_{K}(x, U, g^{-1})=\{k\in K: kg^{-1}(x)\in U\}$ is  right thick in $K$. By Lemma \ref{syn+thick},  $gN_{K}(x, U,g^{-1})g^{-1}\cap K$ is also  right thick in $K$.  Thus  $$gN_{K}(x, U,g^{-1})g^{-1}\cap K \cap g A\neq \emptyset.$$
By noting that $N_{K}(x, U,g^{-1})g^{-1}\subset N_{G}(x, U)$, we have $N_{G}(x, U)\cap A\neq \emptyset$.
\end{proof}

\subsection{Dense regular minimal points}

Let $(X,G)$ be a $G$-system. A point $x\in X$ is {\it regularly recurrent} if for every open neighborhood $U$ of $x$ there is a subgroup $H$ of $G$ with finite index such that $H\subset N_{G}(x, U)$. In particular, a regularly recurrent point is a minimal point.

\medskip
Let $G$ be a discrete group and $(\Gamma_{n})_{n\geq 0}$  a decreasing sequence of subgroups with finite indices in $G$. Let $\pi_{n}: G/\Gamma_{n+1}\rightarrow G/\Gamma_{n}$ be the natural map between coset spaces induced by the inclusion $\Gamma_{n+1}\subset \Gamma_{n}$.  Consider the inverse limit
\[ \overset{\leftarrow}{G}=\underset{\longleftarrow}{\lim}(G/\Gamma_n, \pi_n). \]
Endow every $G/\Gamma_{n}$ with the discrete topology and the product space $\prod_{n\geq 0} G/\Gamma_{n}$ with the product topology. Then  $\overset{\leftarrow}{G}$ is a compact totally disconnected metrizable space, it is a Cantor set if $G/\bigcap_{n} \Gamma_{n}$ is infinite and  finite if $G/\bigcap_{n} \Gamma_{n}$ is finite.

\medskip
Now there is a canonical continuous action of $G$ on $\overset{\leftarrow}{G}$ by left translations, namely for each ${\bf g}=(g_n\Gamma_n)_{n} \in \overset{\leftarrow}{G}$ and $f\in G$,
\[ f\cdot {\bf g}=(fg_n\Gamma_n)_{n}.\]
We denote this $G$-system by $(\overset{\leftarrow}{G}, G)$ and call it a $G$-{\it subodometer} of $G$. If each $\Gamma_n$ is a normal subgroup, then it is called a  $G$-{\it odometer} of $G$.  It is easy to see that $(\overset{\leftarrow}{G}, G)$ is a minimal equicontinuous system and every point is regularly recurrent.

 \medskip
Let $\pi: (X, G)\rightarrow (Y,G)$ be a factor map between two $G$-systems. We say $\pi$ is {\it an almost one-to-one extension}  if  the set $\{x\in X: \pi^{-1}\pi(x)=\{x\}\}$ is dense in $X$.

\medskip
The following result is well known for integer group $\mathbb{Z}$ and it also holds for general groups.
\begin{prop}\cite[Theorem 2]{CP}\label{odometer}
Let $G$ be a discrete group. A minimal $G$-system $(X,G)$ admits a regular minimal point if and only if it is an almost one-to-one extension of a $G$-subodometer.
\end{prop}

The following theorem is an immediate corollary of Theorem \ref{dense small periodic set} and Proposition \ref{odometer}.
\begin{thm}
Let $G$ be a discrete group and $(X,G)$ be a $G$-system. If $(X,G)$ is totally transitive and has dense regular minimal points, then $(X,G)\in \mathcal{M}(G)^{\perp}$.
\end{thm}

In \cite{HY}, the authors construct a transitive $\mathbb{Z}$-subshift that has no periodic points and is disjoint will all minimal $\mathbb{Z}$-systems.

\subsection{Dense distal points}
\begin{lem}\cite[Theorem 6.15]{Aus}\label{distal point}
Let $(X,G)$ and $(Y,G)$ be  $G$-systems. Then for any distal point $x\in X$ and minimal point $y\in Y$, $(x,y)$ is a minimal point in $X\times Y$.
\end{lem}
We remark here that the definition of distal point in \cite{Aus} is a bit different from ours. In \cite{Aus}, a point $x\in X$ is distal if $x$ is the only point in $X$  proximal to itself rather than restricting to the orbit closure of $x$ as we do.  However, Lemma \ref{distal point} still holds by the same proof in \cite[Theorem 6.15]{Aus}, since the point $x'$ constructed in the proof of \cite[Theorem 6.15]{Aus} lies in the orbit closure of $x$.

\begin{thm}
Let $(X,G)$ be an $\infty$-transitive system with a dense set of distal points. Then $(X,G)\in \mathcal{M}(G)^{\perp}$.
\end{thm}
\begin{proof}
Let $\{x_{n}\}_{n=1}^{\infty}$ be a countable dense set of distal points in $X$. Then, by Lemma \ref{proximal cell}, there is an $x\in \bigcap_{n=1}^{\infty}{\rm Prox}_{G}(x_n)$.

\medskip
Let $(Y,G)$ be an arbitrary minimal $G$-system and $J\subset X\times Y$ be a joining. Then there is some $y\in Y$ with $(x,y)\in J$. Fix an $\varepsilon>0$. Then for each $n\geq 1$, the set $\{g\in G: \rho(gx, gx_n)<\varepsilon/2\}$ is right thick in $G$, since $x$ and $x_n$ are proximal. It follows from Lemma \ref{distal point} and the distality of $x_n$ that $(x_n, y)$ is a minimal point in $X\times Y$ and hence $\{g\in G: d(gx_n,x_n)<\varepsilon/2, \rho(gy,y)<\varepsilon\}$ is right syndetic. Thus there is some $g\in G$ such that
\[ \rho(gx_n,x_n)<\varepsilon/2,~~\rho(gx_n, gx)<\varepsilon/2, ~~\rho(gy,y)<\varepsilon.\]
 So $\rho((x_n,y), g(x,y))<\varepsilon$. Since the $\varepsilon$ is arbitrary, we have $(x_n, y)\in \overline{G\cdot (x,y)}$ for each $n\geq 1$. Thus, $X\times\{y\}\subset \overline{G\cdot (x,y)}\subset J$ since $\{x_n\}$ is dense in $X$. By the minimality of $(Y,G)$, we have $J=X\times Y$ and hence $(X,G)$ is disjoint with $(Y,G)$.
\end{proof}

Finally, we show that there exists a countable infinite group
$\Gamma $ such that the set of distal points in the Bernoulli shift $2^{\Gamma}$ is not dense.

\medskip
 Let $\Gamma$ be a Tarski monster group. It is well known that $\Gamma$ is an infinite simple group and  has a property that every nontrivial subgroup of $\Gamma$ is finite.

\medskip
The following are some facts about Tarski monster group and compact groups.
 \begin{lem}\label{Burnside}
 A Tarski monster group cannot be a linear group.
 \end{lem}

 \begin{lem}\label{Lie}\cite[2.20]{MZ}
Let $G$ be a compact group. For every open neighborhood $V$ of the unit ${e_{G}}$, there is a closed normal subgroup $N\subset V$ such that  $G/N$ is a compact linear Lie group.
 \end{lem}

 \begin{lem}\label{Veech}\cite{Veech0}
 Every nontrivial point-distal system admits a nontrivial equicontinuous factor.
 \end{lem}

\begin{thm}\label{equi}
Let $\Gamma$ be a Tarski monster group. Then every minimal equicontinuous $\Gamma$-system is a singleton.
\end{thm}
\begin{proof}
Suppose that $(X,\Gamma)$ is an equicontinuous system. Let $\phi: \Gamma\rightarrow {\rm Homeo}(X)$ denote the action. Then the enveloping semigroup $E(X)$ of  $(X,\Gamma)$  is a compact group and $\phi(\Gamma )$ is dense in $E(X)$ (see \cite[Chapter 3]{Aus}). Let $G=E(X)$.

\medskip
If $ \phi(\Gamma)=\{e_{G}\}$ then we have done. Now assume that $\phi(\Gamma)\neq \{e_{G}\}$. Since $\Gamma$ is a simple group, we have that $\phi(\Gamma)=\Gamma$.
   By Lemma \ref{Burnside} , $G$ is not a compact linear Lie group. It follows from Lemma \ref{Lie} that there is a closed normal subgroup $N\subsetneq G$ such that $G/N$ is a compact Lie group.  But then $N\cap \phi(\Gamma)$ is a proper subgroup of $\phi(\Gamma)$ and hence $N\cap \phi(\Gamma)$  is a finite group. Note that $N\cap \phi(\Gamma)$ is dense in $N$. Thus $N$ is finite. Then we conclude that $G$ is a compact Lie group, which is a contradiction.
   \end{proof}

\begin{cor}\label{tarskinotdense}
Let $\Gamma$ be a Tarski monster group and $(X, \Gamma)$ be a $\Gamma$-system. If $x\in X$ is a distal point then $x$ is fixed by $\Gamma$.  In particular, the set of distal points in the Bernoulli shift $2^{\Gamma}$ is not dense.
\end{cor}
\begin{proof}
Suppose that $x\in X$ is a distal point. Set $Y=\overline{\Gamma x}$. Considering the system $(Y,\Gamma)$, it follows from Lemma \ref{Veech} and Theorem \ref{equi} that $Y$ is a singleton. Hence $x$ is fixed by $\Gamma$.

\medskip
Note that the only fixed points in $2^{\Gamma}$ are ${\bf 1}_{\emptyset}$ and ${\bf 1}_{\Gamma}$. Then the second assertion holds.
\end{proof}

\begin{rem}
Generally, a group admitting a nontrivially minimally equicontinuous action should have a nontrivial homomorphism to a compact group. Such groups cannot be minimally almost periodic 
ones, and  the Bohr compactification of which are trivial. Thus every point-distal system of a minimally almost periodic group is a singleton.
\end{rem}

\section{ $T$-blueprint and proof of Theorem B}
The technics and idea of this section are borrowed from \cite[Chapter 5]{GJS}. In order to apply their results in our setting, we need to generalize their blueprint to $T$-blueprint. The proofs of lemmas in this section are essentially similar to that in \cite[Chapter 5]{GJS} and some proofs are exactly the same. However, we still provide detailed proofs here not only due to completeness but also due to many different notions and conditions.

\subsection{$T$-blueprint and its properties}

\begin{defn}
Let $G$ be a group and let $A,B,\Delta\subset G$. We say that the $\Delta$-translates of $A$ are {\it maximally disjoint within} $B$ if the following properties hold:
\begin{itemize}
\item[(1)] $\Delta$-translates of $A$ are disjoint, i.e., for all $\gamma,\psi\in\Delta$, if $\gamma\neq \psi$ then $A\gamma\cap A\psi=\emptyset$;
\item[(2)] for every $g\in G$, if $Ag\subset B$ then there exists $\gamma\in \Delta$ with $Ag\cap A\gamma \neq \emptyset$.
\end{itemize}
\end{defn}

\begin{defn}
Let $G$ be a countable group and $T\subset G$. A  $T$-{\it preblueprint} is a sequence $(F_n, \Delta_n)_{n\in\mathbb{N}}$ of subsets pair in $G$ satisfying the following conditions:
\begin{itemize}
\item[(i)] (disjoint) for every $n\in\mathbb{N}$ and distinct $\gamma,\psi\in \Delta_n$, $F_n\gamma\cap F_n\psi=\emptyset$;
\item[(ii)] ($T$-covered) for every $n\in\mathbb{N}$, $F_n$ is finite, $e_G\in F_n$, and for every $g\in \Delta_n\setminus\{e_G\}$, $F_ng\subset T$;
\item[(iii)] (coherent) for $k\leq n$, $\gamma\in \Delta_n$ and $\psi\in\Delta_k$,
\[F_k\psi\cap F_n\gamma\neq\emptyset \Longleftrightarrow F_k\psi\subset F_n\gamma; \]
\item[(iv)] (uniform) for $k<n$ and $\gamma,\sigma\in\Delta_n$, $(\Delta_k\cap F_n\gamma)\gamma^{-1}=(\Delta_k\cap F_n\sigma)\sigma^{-1}$;
\item[(v)] (growth) for every $n>0$ and $\gamma\in \Delta_n$, there are distinct $\psi_1,\psi_2,\psi_3\in\Delta_{n-1}$ such that $F_{n-1}\psi_{i}\subset F_{n}\gamma$ for each $i=1,2,3$.
\end{itemize}
A $T$-preblueprint $(F_n, \Delta_n)_{n\in\mathbb{N}}$ is
\begin{itemize}
\item[(1)] a $T$-{\it blueprint} if for every $n\in\mathbb{N}$, $\Delta_n$ is right syndetic;
\item[(2)] {\it maximally disjoint} if the $\Delta_n$-translates of $F_n$ are maximally disjoint within $T$;
\item[(3)] {\it centered} if $e_{G}\in \Delta_n$ for every $n\in\mathbb{N}$;
\item[(4)] {\it directed} if for every $k\in\mathbb{N}$ and $\psi_1, \psi_2\in\Delta_k$ there is $n>k$ and $\gamma\in \Delta_n$ with $F_k\psi_1\cup F_k\psi_2\subset F_n\gamma$.
\end{itemize}
\end{defn}

 By the item (iv) above, we denote the set $(\Delta_k\cap F_n\gamma)\gamma^{-1}$ by $D^{n}_{k}$, since it is independent on the choice of $\gamma$.

\begin{lem}\label{property of blueprint}
Let $G$ be a countable group and  $T\subset G$. Suppose that $(F_n, \Delta_n)_{n\in\mathbb{N}}$ is a  $T$-preblueprint.
\begin{itemize}
\item[(a)] If $(F_n, \Delta_n)_{n\in\mathbb{N}}$ is maximally disjoint and $T$ is right thickly syndetic then $(F_n, \Delta_n)_{n\in\mathbb{N}}$ is a $T$-blueprint.
\item[(b)] If $(F_n, \Delta_n)_{n\in\mathbb{N}}$ is directed, then for any $r(1), r(2),\ldots, r(m)$ and $\psi_1\in \Delta_{r(1)},\ldots, \psi_m\in\Delta_{r(m)}$, there is $n\in \mathbb{N}$ and $\gamma \in\Delta_{n}$ such that for every $1\leq i\leq m$, $F_{r(i)}\psi_{i}\subset F_n\gamma$.
\item[(c)] If $t\leq k<n$, $\sigma\in\Delta_n, A\subset G$ is finite and satisfies $F_tF_kA\subset T$, $\Delta_t\cap F_t^{-1}F_t F_k A\subset D_t^n\sigma $, and the $\Delta_t$-translates of $F_t$ are maximally disjoint within $T$, then for all $\gamma\in \Delta_n$,
\[ \Delta_k\cap A\sigma^{-1}\gamma=(\Delta_{k}\cap A)\sigma^{-1}\gamma\subset  D_{k}^{n}\gamma.\]
\item[(d)] If $(F_n, \Delta_n)_{n\in\mathbb{N}}$ is is a directed $T$-blueprint, the $\Delta_0$-translates of $F_0$ are maximally disjoint within $T$ and $T$ is right thickly syndetic, then it is minimal in the following sense: for every finite subset $A\subset G$  there is $h\in G$ and a finite subset $F\subset G$ such that for any
\[\forall g\in G~~\exists f\in F,~~(\Delta_0\cap Ah)fg=\Delta_{0}\cap Ahfg.\]
\end{itemize}
\end{lem}
\begin{proof}
(a)  By the definition, we need to show that each $\Delta_n$ is right syndetic. Fix $n\in\mathbb{N}$. Since $T$ is right thickly syndetic, there is a subset $Q_n\subset G$ and a finite subset $A_n\subset G$ satisfying
 \[ F_nQ_n\subset T~~~~ {\rm and }~~~~ A_nQ_n=G.\]
Given $g\in G$, we write $g=aq$ with $a\in A_n$ and $q\in Q_n$. Thus $q=a^{-1}g$ and hence $F_na^{-1}g\subset T$. Now that $(F_n, \Delta_n)_{n\in\mathbb{N}}$ is maximally disjoint, we have $F_na^{-1}g\cap F_n\Delta_n\neq \emptyset$ and hence
\[g\in aF_n^{-1}F_n\Delta_n\subset A_nF_n^{-1}F_n\Delta_n.\]
Since $g$ is arbitrary, we have $G=A_nF_n^{-1}F_n\Delta_n$ and hence $\Delta_n$ is right syndetic.

\medskip
(b) It suffices to show that the clause holds for the maximal elements with respect to inclusion among $\{F_{r(1)}\psi_{1}, \ldots, F_{r(m)}\psi_{m}\}$.  By the coherent property, distinct maximal elements of this collection are disjoint. Thus we may assume that $F_{r(1)}\psi_{1}, \ldots, F_{r(m)}\psi_{m}$ are mutually disjoint and we may further assume that $r(1)\leq r(2)\leq\cdots\leq r(m)$. For each $i>1$, pick $\lambda_i\in D_{r(1)}^{r(i)}$. Then for each $i>1$, we have $\lambda_i\psi_i\in \Delta_{r(1)}$. Now for each $i>1$, choose $n(i)$ and $\sigma(i)\in \Delta_{n(i)}$ with
\[ F_{r(1)}\psi_1\cup F_{r(1)}\lambda_i\psi_i\subset F_{n(i)}\sigma_{i}.\]
Thus for $i,j>1$, we have
\[ F_{r(1)}\psi_1 \subset    F_{n(i)}\sigma_{i}\cap  F_{n(j)}\sigma_{j}.\]
In particular, $F_{n(i)}\sigma_{i}$ and $F_{n(j)}\sigma_{j}$ have nonempty intersection. By the coherent property, it must be that one of $F_{n(i)}\sigma_{i}$ and $F_{n(j)}\sigma_{j}$ contains the other one. Then there is some $n\in\{n(2), \ldots, n(m)\}$ and $\sigma\in \Delta_{n}$ such that $F_{n(i)}\sigma \subset F_n\sigma$ for each $i>1$. Thus
\[   F_{r(1)}\psi_1\cup F_{r(1)}\lambda_2\psi_2\cup\cdots\cup F_{r(1)}\lambda_m\psi_m\subset F_n\sigma.\]
Since $F_{r(1)}\lambda_i\psi_i\subset F_{r(1)}\psi_i$ for each $i>1$, we have $F_{r(1)}\psi_i\cap F_n\sigma\neq \emptyset$. However, none of $F_{r(i)}\psi_i$'s can contain $F_n\sigma$ since they are pairwise disjoint. Therefore, $F_n\sigma$ contain each $F_{r(i)}\psi_i$ by coherent property.

\medskip
(c)  Fix $\gamma\in \Delta_{n}$. Then we have
\[ D_{t}^{k}(\Delta_{k}\cap A)\subset D_{t}^{k}\Delta_{k}\cap D_{t}^{k}A
\subset \Delta_{t}\cap F^{-1}_{t}F_tF_k A\subset D_{t}^{n}\sigma\subset F_n\sigma. \]
Thus we have $F_{k}\psi\cap F_n\sigma\neq \emptyset$ for any $\psi\in \Delta_{k}\cap A$. By the coherent property, we have $\Delta_k\cap A\subset D_{k}^n \sigma$. So
\[ (\Delta_k\cap A)\sigma^{-1}\gamma\subset D_{k}^{n}\sigma\sigma^{-1}\gamma=D_{k}^{n}\gamma\subset D_{k}^{n}\Delta_{n}\subset \Delta_{k}.\]
It is trivial that $(\Delta_k\cap A)\sigma^{-1}\gamma\subset A \sigma^{-1}\gamma$. Therefore,
\[ (\Delta_{k}\cap A)\sigma^{-1}\gamma\subset \Delta_k\cap A\sigma^{-1}\gamma.\]

\medskip
To show the reverse inclusion, take $\lambda\in \Delta_k\cap A\sigma^{-1}\gamma$ and fix any $\tau\in D_{t}^{k}$. Then we have $\tau\lambda\gamma^{-1}\sigma\in F_{k}A$ and hence
\[ F_t\tau\lambda\gamma^{-1}\sigma \subset F_tF_k A.\]
{\bf Claim}. $\Delta_{t}\cap F^{-1}_tF_tF_k A$-translates of $F_t$ are maximally disjoint within $F_tF_kA$.

Indeed, for any $g\in G$, if $F_tg\subset F_tF_kA\subset T$ then there is some $h\in \Delta_{t}$ such that $F_tg\cap F_th\neq\emptyset$. Thus $h\in F_t^{-1}F_tg\subset F_t^{-1}F_tF_kA $. So $h\in \Delta_t\cap F_t^{-1}F_tF_kA$ and hence the claim holds.

Now the claim implies that there is some $\psi\in \Delta_t\cap F_t^{-1}F_tF_kA$ such that
$ F_t\psi\cap F_t\tau\lambda\gamma^{-1}\sigma\neq\emptyset$ and hence
\begin{equation}\label{eq 0}
F_t\psi\sigma^{-1}\gamma\cap F_t\tau\lambda\neq\emptyset.
\end{equation}
By noting that $\psi\in D_{t}^{n}\sigma$, we have $\psi\sigma^{-1}\gamma\in D_{t}^{n} \gamma=\Delta_t\cap F_n\gamma$. Thus $\psi\sigma^{-1}\gamma\in\Delta_t$. Also $\tau\lambda\in D_{t}^{k}\lambda=\Delta_t\cap F_k\lambda$ and hence $\tau\lambda\in\Delta_t$. This together with (\ref{eq 0}) imply that $\psi\sigma^{-1}\gamma=\tau\lambda$.
Note that we also have $\psi\sigma^{-1}\gamma\in F_n\gamma$ and $\tau\lambda\in F_k\lambda$.  Thus $F_n\gamma\cap F_k\lambda\neq\emptyset$. By the coherency property, we have $F_k\lambda\subset F_n\gamma$ and hence $\lambda\in D_{k}^{n}\gamma$. Thus $\lambda\gamma^{-1}\sigma\in D_{k}^{n}\sigma\subset \Delta_k$. So $\lambda\gamma^{-1}\sigma\in \Delta_k\cap A$ and hence $\lambda\in(\Delta_k\cap A)\sigma^{-1}\gamma$. This shows the reverse inclusion.

\medskip
(d) Let $A$ be a finite subset of $G$. Since $T$ is right thick, there is some $h\in G$ such that $F_0F_0Ah\subset T$. Let $C=\Delta_0\cap F_0^{-1}F_0F_0Ah$. Now (b) implies that there is $n\in\mathbb{N}$ and $\sigma\in\Delta_n$ such that $F_0C\subset F_n\sigma$. Particularly, $C\subset D_{0}^n\sigma$. Since $(F_n, \Delta_n)$ is a $T$-blueprint, $\Delta_n$ is right syndetic and hence there is a finite subset $B\subset G$ with $B\Delta_n=G$.

We claim that $F=\sigma^{-1}B^{-1}$ satisfies the requirements. Fix $g\in G$. Since $B\Delta_n=G$, there is $b\in B$ and $\gamma\in \Delta_n$ with $g=b\gamma$ and hence $b^{-1}g=\gamma\in\Delta_n$. Applying (c) to $Ah$ and $\sigma, \gamma$, we have
\[ (\Delta_0\cap Ah)\sigma^{-1}\gamma=\Delta_0\cap Ah\sigma^{-1}\gamma.\]
By setting $f=\sigma^{-1}b^{-1}\in F$, we have $fg=\sigma^{-1}b^{-1}g=\sigma^{-1}\gamma$. Thus
\[ (\Delta_0\cap Ah)fg=\Delta_0\cap Ahfg.\]
 \end{proof}

The proof of the following theorem is given in the next subsection.
\begin{thm}\label{existence of blueprint}
Let $G$ be a countably infinite group and $T\subset G$ be right thickly syndetic. Then for any finite subset $A\subset G$ with $e_G\in A$,   there is a maximally disjoint, centered, directed $T$-blueprint $(F_n, \Delta_n)_{n\in\mathbb{N}}$ satisfying $F_0=A$.
\end{thm}

\begin{lem}\label{char of min}
Let $G$ be a countable group and $x\in 2^{G}$. Then $x$ is an almost periodic point in $2^{G}$ if and only if for every finite subset $A\subset G$ there exist $h\in G$ and  a finite subset $F\subset G$ such that
\begin{equation}\label{eq 3}
 \forall g\in G~~\exists f\in F~~\forall a\in A,~~~x(ahfg)=x(ah).
 \end{equation}
\end{lem}
\begin{proof}
For  $y\in 2^G$ and $B\subset G$, let
\[ B(y):=\{ z\in 2^G: \forall g\in B, ~~z(g)=y(g)\}.\]

Suppose that $x$ is an almost periodic point. Then $N(x, A(x))=\{g\in G: gx\in A(x)\}$ is right syndetic. Thus there is a finite subset $F\subset G$ such that $F^{-1}N(x, A(x))=G$. So for any $g\in G$, there is some $f\in F$ so that $fg\in N(x, A(x))$, i.e., $x(a)=fgx(a)=x(afg)$ for any $a\in A$.

\medskip
Suppose that for any finite set $A\subset G$, there is $h\in G$ and finite $F\subset G$ so that (\ref{eq 3}) holds. For $g\in G$ and $f\in F$, if $x(ahfg)=x(ah)$ for any $a\in A$, then $hx(ahfgh^{-1})=hx(a)$. So $hfgh^{-1}\in N(hx, A(hx))$. Thus we have $ G=F^{-1}h^{-1}N(hx, A(hx))h$.
It follows that $$G=Gh^{-1}=F^{-1}h^{-1}N(hx, A(hx)).$$ So $N(hx, A(hx))$ is right syndetic. Since $A$ is arbitrary, we conclude that $hx$ is almost periodic and hence so is $x$.
\end{proof}

Now we show Theorem B.
\begin{proof}[Proof of Theorem B]
We may assume that $G$ is countably infinite. Now choose a finite nonempty subset $A\subset T$ and let $(F_n, \Delta_n)_{n\in\mathbb{N}}$ be a directed maximally disjoint $T$-blueprint with $A\subset F_0$.

\medskip
 Let $S=A\Delta_{0}$. By the definition of $T$-blueprint and noting that $A\subset T$, we have $S\subset T$. It remains to show that ${\bf 1}_{S}$ is an almost periodic point in $2^{G}$. Let $B$ be a finite subset of $G$.  According to Lemma \ref{property of blueprint} (d), there is $h\in G$ and a finite subset $F\subset G$ such that for any $g\in G$, there exists $f\in F$ such that $(\Delta_0\cap A^{-1}Bh)fg=\Delta_0\cap A^{-1}Bhfg$. This is equivalent to say that for any $g\in G$ there is $f\in F$ such that
 \begin{equation}\label{eq 1}
  \forall c\in A^{-1}B~~~~(chfg\in \Delta_0~~\Longleftrightarrow~~ch\in\Delta_0).
  \end{equation}

 Now let $g\in G$ be arbitrary and let $f\in F$ satisfy (\ref{eq 1}). We claim that
 \begin{equation}\label{eq 2}
 \forall b\in B~~~~(bh\in S~~\Longleftrightarrow~~bhfg\in S).
 \end{equation}
Let  $b\in B$.  If $bhfg=a\gamma\in S$ for some $a\in A$ and $\gamma\in\Delta_0$, then $a^{-1}bhfg=\gamma\in \Delta_0$. Thus, by (\ref{eq 1}),  $a^{-1}bh\in \Delta_0$  and hence $bh\in a\Delta_0\subset S$. Similarly, if $bh\in a\Delta_0$ for some $a\in A$, then $a^{-1}bh\in \Delta_0$ and hence $a^{-1}bhfg\in \Delta_0$ by (\ref{eq 1}). Thus $bhfg\in a\Delta_0\subset S$. This shows (\ref{eq 2}). By Lemma \ref{char of min},  (\ref{eq 2}) implies the almost periodicity of ${\bf 1}_{S}$.
\end{proof}

 Now we are going to prove Proposition \ref{m-set}.
 \begin{proof}[Proof of Proposition \ref{m-set}]
Let $T$ be a right thickly syndetic subset of $G$. By Theorem B, 
there is a subset $S\subset T$ such that ${\bf 1}_{S}$ is a minimal point in $2^{G}$. Let $U=\{x\in 2^{G}: x(e_G)=1\}$. Then
\[ N({\bf 1}_{S}, U)=\{ g\in G: x(g)=1\}= S.\]
Thus $S$ is an $m$-set and contained in $T$.
\end{proof}

\subsection{Existence of $T$-blueprint}
\begin{lem}\label{construction of preblueprint}
Let $G$ be a countably infinite group and $T\subset G$. Suppose that $(F_n)_{n\in\mathbb{N}}$ is a sequence of finite subsets of $G$ and $\{\delta_{k}^{n}: n,k\in \mathbb{N}, k<n\}$ is a collection of finite subsets of $G$ satisfying
\begin{itemize}
\item[(i)] $e_{G}\in \delta_{n-1}^{n}$ for each $n\geq 1$;
\item[(ii)] $|\delta_{n-1}^{n}|\geq 3$ for each $n\geq 1$;
\item[(iii)] for all $n,k\in\mathbb{N}$ with $k<n$, the $\delta_{k}^{n}$-translates of $F_k$ are disjoint;
\item[(iv)] for all $n,k\in\mathbb{N}$ with $k<n$, and every choice of $k<i_1<i_2<\cdots<i_{t}<n$,
\[F_{k} \left(((\delta_{k}^{n}\cup(\delta_{k}^{i_1}\delta_{i_1}^{i_2}\cdots \delta_{i_{t-1}}^{i_t}\delta_{i_t}^{n})) \setminus \{e_G\}\right)\subset T;\]
\item[(v)] $F_m\delta_{m}^{n}\cap F_k\delta_{k}^{n}=\emptyset$ for all $m\neq k<n$;
\item[(vi)] $F_n=\bigcup_{0\leq k<n}F_k\delta_{k}^{n}$ for all $n\geq 1$.
\end{itemize}
Then there is a sequence $(\Delta_n)_{n\in\mathbb{N}}$ of subsets of $G$ with $\delta_{k}^{n}\subset \Delta_{k}$ for every $n,k\in \mathbb{N}$ with $k<n$ and such that $(F_n, \Delta_n)_{n\in\mathbb{N}}$ is a centered and directed $T$-preblueprint.
\end{lem}

\begin{proof}
For $k,n\in\mathbb{N}$ with $k<n$, define $D_{k}^{k}=\{e_G\}, D_{k}^{k+1}=\delta_{k}^{k+1}$, and inductively define
\begin{eqnarray*}
D_{k}^{n}&=&\delta_{k}^{n}\cup D_{k}^{k+1}\delta_{k+1}^{n}\cup\cdots\cup D_{k}^{n-2}\delta_{n-2}^{n}\cup D_{k}^{n-1}\delta_{n-1}^{n}=
\bigcup_{k\leq m<n}D_{k}^{m}\delta_{m}^{n}\\
&=&\delta_{k}^{n}\cup\bigcup_{t=1}^{n-k}  \left\{\delta_{k}^{i_1}\delta_{i_1}^{i_2}\delta_{i_2}^{i_3}\cdots\delta_{i_{t-1}}^{i_t}\delta_{i_{t}}^{n}:~k<i_1<i_2<\cdots<i_t<n\right\} .
\end{eqnarray*}

\noindent{\bf Claim 1}. $F_kD_{k}^{n}\subset F_n$ for all $k,n\in\mathbb{N}$ with $k\leq n$.

\medskip
It is clear that $F_kD_{k}^{k}=F_k$. If we assume that $F_kD_{k}^{m}\subset F_m$ for all $k\leq m<n$, then
\[ F_kD_{k}^{n}=\bigcup_{k\leq m<n} F_k D_{k}^{m}\delta_{m}^{n} \subset \bigcup_{k\leq m<n} F_{m}\delta_{m}^{n}\subset F_n. \]
Now we can conclude the claim by induction.

\medskip
\noindent{\bf Claim 2}. $D_{k}^{m}D_{m}^{n}\subset D_{k}^{n}$ for all $k,m,n\in\mathbb{N}$ with $k\leq m\leq n$.

\medskip
When $m=n$, we have $D_{k}^{m}D_{m}^{n}=D_{k}^{n}D_{n}^{n}=D_{k}^{n}$. If we assume that $D_{k}^{r}D_{r}^{n}\subset D_{k}^{n}$ for all $m\leq r<n$, then
\[ D_{k}^{m}D_{m}^{n}=\bigcup_{m\leq r<n} D_{k}^{m}D_{m}^{r}\delta_{r}^{n}\subset \bigcup_{m\leq r<n} D_{k}^{r}\delta_{r}^{n}\subset \bigcup_{k\leq r<n} D_{k}^{r}\delta_{r}^{n}\subset D_{k}^{n}.\]
Now the claim follows from induction.

\medskip
\noindent{\bf Claim 3}. The $D_{k}^{n}$-translates of $F_k$ are disjoint for all $k,n\in\mathbb{N}$ with $k\leq n$.

\medskip
When $n=k$ or $n=k+1$, it is clear that the $D_{k}^{n}$-translates of $F_k$ are disjoint by (iii). Now we assume that the $D_{k}^{m}$-translates of $F_k$ are disjoint for each $m$ with $k\leq m<n$  and are going to $D_{k}^{n}$-translates of $F_k$ are disjoint.  Let $f,g\in D_{k}^{n}$ be two distinct elements.  We need to show $F_kf\cap F_kg=\emptyset$. Recall that $D_{k}^{n}=
\bigcup_{k\leq m<n}D_{k}^{m}\delta_{m}^{n}$.
\begin{itemize}
\item[(1)] If $f\in D_{k}^{r}\delta_{r}^{n}$ and $g\in D_{k}^{s}\delta_{s}^{n}$ for some $k\leq r<s<n$, then by Claim 1, $F_kf\subset F_r\delta_{r}^{n}$ and $F_kg\subset F_s\delta_{s}^{n}$. It follows from  (v) that $F_kf\cap F_kg=\emptyset$.
\item[(2)] If $f\in  D_{k}^{m}\gamma$ and $g\in  D_{k}^{m}\psi$ for some $k\leq m<n$ and distinct $\gamma, \lambda\in \delta_{m}^{n}$, then $F_kf\subset F_m\gamma$ and $F_kg\subset F_m\psi$ by Claim 1. Thus  it  also follows from  (v) that $F_kf\cap F_kg=\emptyset$.
\item[(3)] If $f,g\in D_{k}^{m}\gamma$ for some $k\leq m<n$ and $\gamma\in \delta_{m}^{n}$, then $F_kf\cap F_kg=\emptyset$ by our assumption that $D_{k}^{m}$-translates of $F_k$ are disjoint.
\end{itemize}

\medskip
By (iv) and the definition of $D_{k}^{n}$, we conclude the following claim.

\noindent{\bf Claim 4}. For all $k,n\in\mathbb{N}$ with $k<n$ and any $g\in D_{k}^{n}\setminus\{e_{G}\}$, $F_kg\subset T$.

\bigskip
For $n,k\in\mathbb{N}$ with $k\leq n$, we have $D_{k}^{n}\subset D_{k}^{n+1}$ since $D_{k}^{n}\delta_{n}^{n+1}$ and $e_{G}\in \delta_{n}^{n+1}$. For each $k\in \mathbb{N}$, we define
\[\Delta_{k}=\bigcup_{n\geq k} D_{k}^{n}.\]
Next we will verify that $(F_n, \Delta_n)_{n\in\mathbb{N}}$ is a $T$-preblueprint.

\medskip
(Disjoint) Let $n\in\mathbb{N}$ and $\gamma,\psi\in \Delta_{n}$. Then there is some $m>n$ such that $\gamma,\psi\in D_{n}^{m}$ by the monotonicity of $D_{n}^{m}$'s. It now follows from the Claim 3 that $F_{n}\gamma\cap F_{n}\psi=\emptyset$.

\medskip
($T$-covered) This follows from the Claim 4.

\medskip
(Coherent) Suppose that $F_{k}\psi\cap F_{n}\gamma\neq\emptyset$, for some $k<n$ and $\psi\in\Delta_{k}, \gamma\in \Delta_{n}$. Let $m\geq n$ be large enough so that $\psi\in D_{k}^{m}$ and $\gamma\in D_{n}^{m}$. To show $F_{k}\psi\subset F_{n}\gamma$, it suffices to show  that $\psi\in D_{k}^{n}\gamma$, since $F_{k}D_{k}^{n}\subset F_{n}$.

We show $\psi\in D_{k}^{n}\gamma$ by induction on $m$. It is trivial when $m=n$ since then $\gamma=e_G$ and hence $\psi\in D_{k}^{n}=D_{k}^{n}\gamma$. Suppose $m>n$ and the claim holds for any $i$ with $n\leq i<m$. According to the definition of $D_{k}^{m}$ and $D_{n}^{m}$, there are $s, t$ with $k\leq s<m$ and $n\leq t<m$ such that $\psi\in D_{k}^{s}\delta_{s}^{m}$ and $\gamma\in D_{n}^{t}\delta_{t}^{m}$. If $s\neq t$ then (v) implies
\[ F_k\psi\cap F_n\gamma\subset F_{k}D_{k}^{s}\delta_{s}^{m}\cap F_{n}D_{n}^{t}\delta_{t}^{m}\subset F_{s}\delta_{s}^{m}\cap F_{t}\delta_{t}^{m}=\emptyset.\]
Thus we must have $s=t$. Let $\lambda,\sigma\in \delta_{s}^{m}$ be such that $\psi\in D_{k}^{s}\lambda$ and $\gamma\in D_{n}^{s}\sigma$. If $\lambda \neq \sigma$, then (iii) implies
\[ F_k\psi\cap F_n\gamma\subset F_{k}D_{k}^{s}\lambda\cap F_{n}D_{n}^{s}\sigma\subset F_{s}\lambda\cap F_{s}\sigma=\emptyset.\]
Thus $\lambda=\sigma$ and hence $ \psi\lambda^{-1}\in D_{k}^{s}\subset  \Delta_{k}$, $\gamma\lambda^{-1}\in D_{n}^{s}\subset\Delta_n$, and
\[ F_{k}\psi\lambda^{-1}\cap F_{n}\gamma\lambda^{-1}\neq \emptyset.\]
By the induction hypothesis, we have $\psi\lambda^{-1}\in D_{k}^{n}\gamma\lambda^{-1}$ and hence $\psi\in D_{k}^{n}\gamma$.

\medskip
(Uniform) It suffices to show $\Delta_{k}\cap F_{n}\gamma=D_{k}^{n}\gamma$ for all $k<n$ and $\gamma\in \Delta_{n}$.  Let $m>n$ be such that $\gamma\in D_{n}^{m}$. Then it follows from Claim 2 that $D_{k}^{n}\gamma\subset D_{k}^{n}D_{n}^{m}\subset D_{k}^{m}\subset \Delta_{k}$. By Claim 1, $F_{k}D_{k}^{n}\gamma\subset F_n\gamma$. Thus
\[ D_{k}^{n}\gamma\subset \Delta_{k}\cap F_{n}\gamma.\]
For the converse inclusion, if $\psi\in \Delta_{k}\cap F_{n}\gamma$ then $F_{k}\psi\cap F_{n}\gamma\neq \emptyset$ since $e_{G}\in F_{k}$. Then it follows from the proof of the coherency that $\psi\in D_{k}^{n}$.

\medskip
(Growth) $|D_{n-1}^{n}|=|\delta_{n-1}^{n}|\geq 3$.

\medskip
(Centered) $e_G\in \delta_{n}^{n+1}\subset\Delta_n$.

\medskip
(Directed) Let $n,k\in\mathbb{N}$ and $\gamma\in \Delta_n, \psi\in\Delta_k$. Then there is some large $m$ with $\gamma\in D_{n}^{m}$ and $\psi\in D_{k}^{m}$. Thus $F_{n}\gamma, F_{k}\psi\subset F_{m}e_{G}$.

\end{proof}

\begin{defn}
Let $G$ be a countably infinite group and $T\subset G$. A {\it growth sequence with respect to $T$ } is a sequence $(H_n)_{n\in\mathbb{N}}$ of finite subsets of $G$ satisfying
\begin{itemize}
\item[(i)] $e_{G}\in H_0$;
\item[(ii)] $\bigcup_{n\in\mathbb{N}}H_{n}=G$;
\item[(iii)] $(H_{n-1}H_{n-1}^{-1})\cdots (H_{1}H_{1}^{-1})(H_0H_0^{-1})H_{n-1}\subset H_{n}$ for every $n\geq 1$;
\item[(iv)] for each $n\geq 1$, if $\Delta\subset H_{n}$ satisfies $H_{n-1}\Delta\subset H_n\cap T$  and has the property that $H_{n-1}g\cap H_{n-1}\Delta\neq \emptyset$ whenever $H_{n-1}g\subset H_{n}\cap T$, then $|\Delta|\geq 3$.
\end{itemize}
\end{defn}

\begin{lem}\label{existence of growth sequence}
Let $G$ be a countably infinite group and $A\subset G$ be finite with $e_G\in A$. Then, for any  right  syndetic set $T\subset G$,  there is a growth sequence $(H_n)_{n\in\mathbb{N}}$ with respect to $T$ such that $H_0=A$.
\end{lem}
\begin{proof}
Let $g_0=e_G, g_1, g_2,\ldots$ be an enumeration of $G$. Let $H_0=A$.  Now suppose that we have defined $H_0, \ldots H_{n-1}$ for some $n\geq 1$. Let $k=3|H_{n-1}|(|H_{n-1}|^2+1)$. By Lemma \ref{dis trans}, there are at least $3|H_{n-1}|$ disjoint translates of $H_{n-1}$ among $$H_{n-1}g_0, H_{n-1}g_1,\ldots, H_{n-1}g_{k-1}.$$
Since $T$ is right thickly, there is some $h\in G$ such that $\bigcup_{i=0}^{k-1} H_{n-1}g_ih\subset T$. Define
\[ H_{n}= \{g_0,\ldots,g_n\}\cup (H_{n-1}H_{n-1}^{-1})\cdots (H_{1}H_{1}^{-1})(H_0H_0^{-1})H_{n-1}\cup \bigcup_{i=0}^{k-1} H_{n-1}g_ih.\]
Clearly, $(H_n)_{n\in \mathbb{N}}$ satisfies (i),(ii) and (iii)  in the definition of growth sequence.
Note that for each $n\geq 1$, $H_n\cap T$ contains at least $3|H_{n-1}|$ disjoint translates of $H_{n-1}$. Thus, if $\Delta\subset H_{n}$ satisfies that $H_{n-1}\Delta\subset H_{n}$ and has the property that $H_{n-1}g\cap H_{n-1}\Delta\neq \emptyset$ whenever $H_{n-1}g\subset H_n$, then $|\Delta|\geq 3$ since $H_{n}\Delta$ can intersect at most $ |H_{n-1}||\Delta|$ translates of $H_{n-1}$. This establishes (iv) and hence $(H_{n})_{n\in\mathbb{N}}$ satisfies the requirements.
\end{proof}

\begin{thm}\label{existence of blueprint under growth sequence}
Let $G$ be a countably infinite group and $T\subset G$ be a right thickly syndetic subset. Let $(H_n)_{n\in\mathbb{N}}$ be a growth sequence with respect to $T$. Then there is a maximally disjoint, centered, directed $T$-blueprint $(F_n, \Delta_n)_{n\in\mathbb{N}}$ satisfying
\begin{itemize}
\item[(a)] $F_0=H_0$;
\item[(b)] $F_n\subset H_n$ for each $n\geq 1$;
\item[(c)] for each $n\geq 1$, the $D_{n-1}^{n}$-translates of $F_{n-1}$ are contained and maximally disjoint within $H_n\cap T$;
\item[(d)]  for all $n\geq 1$ and $0\leq k<n$, the $D_{k}^{n}$-translates of $F_{k}$ are  maximally disjoint within $H_{n-1}$.
\end{itemize}
\end{thm}

\begin{proof}
Set $F_0=H_0$ so that (a) holds.  Fix $n\geq 1$ and we assume that we have constructed $F_0, \ldots, F_{n-1}$ with each $F_i\subset H_i$.  Choose $\delta_{n-1}^{n}$ such that
\begin{itemize}
\item $e_{G}\in \delta_{n-1}^{n}$ and
\item the $(\delta_{n-1}^{n}\setminus\{e_G\})$-translates of $F_{n-1}$ are contained and maximally disjoint within $H_n\cap T$.
\end{itemize}
By the properties of growth sequence, we must have $|\delta_{n-1}^{n}|\geq 3$. Now assume that we have defined $\delta_{n-1}^{n}$ through $\delta_{k+1}^{n}$ with $0\leq k<n-1$.

Let \[ T_{k}^{n}=\left\{g\in T: ~~(F_{n-1}F_{n-1}^{-1})\cdots(F_{k+2}F_{k+2}^{-1})(F_{k+1}F_{k+1}^{-1})g\subset H_n\right\},\]
and let $\delta_{k}^{n}$ be such that the $\delta_{k}^{n}$-translates of $F_{k}$ are contained and maximally disjoint within $T_{k}^{n}-\bigcup_{k<m<n}F_{m}\delta_{m}^{n}$.  Finally, define
\[ F_n=\bigcup_{0\leq k<n}F_{k}\delta_{k}^{n}.\]

\noindent{\bf Claim 1}. The $F_{n}$'s and $\delta_{k}^{n}$'s satisfy the conditions of Lemma \ref{construction of preblueprint}.

\medskip
It is clear that the $F_{n}$'s and $\delta_{k}^{n}$'s satisfy all conditions of Lemma \ref{construction of preblueprint} except (iv). To show (iv), fix $n,k\in\mathbb{N}$ with $k<n$, and  $k<i_1<i_2<\cdots<i_{t}<n$. By definition of $\delta_{n-1}^{n}$, we may further assume $k<n-1$. Then, by the definition of  $\delta_{k}^{n}$, we have $F_k\delta_{k}^{n}\subset T_{k}^{n}\subset T$. For $g\in \delta_{k}^{i_1}\delta_{i_1}^{i_2}\cdots \delta_{i_{t-1}}^{i_t}\delta_{i_t}^{n}\setminus\{e_G\}$, write $g=g_0g_1\cdots g_t$ with $g_0\in \delta_{k}^{i_1}, g_1\in \delta_{i_1}^{i_2}, \ldots, g_t\in \delta_{i_t}^{n}$.  If $g_t\neq e_G$, then
\[F_kg=F_{k}g_0g_1\cdots g_t\subset F_{i_1}g_1\cdots g_t\subset\cdots\subset F_{i_{t}}g_t\]
and the definition of $\delta_{i_t}^{n}$ implies $F_{i_t}g_t\subset H_{n}\cap T$. Thus $F_{k}g\subset T$. If $g_{t}=e_G$,  let $r\in\{1,\ldots,t\}$ be the largest one such that $g_{r}\cdots g_{t}\neq e_{G}$. Then we have $g=g_0\cdots g_r$ with $g_r\neq e_G$, and
\[F_kg=F_{k}g_0g_1\cdots g_{r-1} g_r\cdots g_t=F_kg_0\cdots g_r  \subset   F_{i_{r}}g_{r}\]
and we also have $F_kg\subset H_{i_{r}}\cap T\subset T$. This shows the claim.

\medskip
Now Lemma \ref{construction of preblueprint} implies that $(F_n,\Delta_n)_{n\in\mathbb{N}}$ is a centered and directed $T$-preblueprint with $\Delta_n$ as defined in Lemma   \ref{construction of preblueprint}. It is clear that $(F_n,\Delta_n)_{n\in\mathbb{N}}$ satisfies (a), (b) and (c). Next we will show that it also satisfies (d).

\medskip
\noindent{\bf Claim 2}. If $n>k, g\in G, F_{k}g\subset T$, and $F_{k}g\cap F_n\neq \emptyset$, then $F_kg\cap F_m\delta_{m}^{n}\neq\emptyset$ for some $m$ with $k\leq m<n$.

 \medskip
Now we may assume that  $F_kg\cap F_{m}\delta_{m}^{n}=\emptyset$ for all $k<m<n$ and we are going to show $F_{k}g\cap F_{k}\delta_{k}^{n}\neq\emptyset$,  which then implies the claim holds with $m=k$.  Since $F_n=\bigcup_{0\leq t<n}F_t\delta_{t}^{n}$, there is some $0\leq t\leq k$ with $F_kg\cap F_{t}\delta_{t}^{n}\neq\emptyset$. Thus $g\in F_{k}^{-1}F_{t}\delta_{t}^{n}$.
We may further assume that $t<k$, otherwise we are done. Then we have
\[ F_{k}g\subset F_{k}F_{k}^{-1}F_{t}\delta_{t}^{n}\subset (F_{k}F_{k}^{-1})\cdots(F_{t+1}F_{t+1}^{-1})F_{t}\delta_{t}^{n}.\]
So
 \[ (F_{n}F_{n}^{-1})\cdots(F_{k+1}F_{k+1}^{-1})F_{k}g\subset (F_{n}F_{n}^{-1})\cdots(F_{t+1}F_{t+1}^{-1})F_{t}\delta_{t}^{n}.\]
By noting that $F_{t}\delta_{t}^{n}\subset T_{t}^{n}$, the right hand side of the above is contained in $H_{n}$, and hence $F_kg\subset T_{k}^{n}$ since $F_{k}g\subset T$. Therefore,
\[ F_{k}g\subset T_{k}^{n}-\bigcup_{k<m<n}F_{m}\delta_{m}^{n}.\]
  It now follows from the definition of $\delta_{k}^{n}$ that $F_kg\cap F_{k}\delta_{k}^{n}\neq \emptyset$. Thus we have showed the Claim 1.

\medskip
\noindent{\bf Claim 3}.  For all $g\in G$ and $k\leq n$, if $F_kg\cap F_n\neq \emptyset$ then  $F_{k}\cap F_{k}D_{k}^{n}\neq\emptyset$, where $D^{n}_{k}$ is defined as in Lemma \ref{construction of preblueprint}.

\medskip
Fix $k\in\mathbb{N}$. We may assume that  $n>k$ and the claim holds for any $k\leq m<n$. Let $g\in G$ be such that $F_kg\cap F_n\neq \emptyset$.  Then, by Claim 2, there is some $m$ with $k\leq m<n$ such that $F_kg\cap F_m\delta_{m}^{n}\neq \emptyset$. Let $\gamma\in \delta_{m}^{n}$ be such that $F_kg\cap F_m\gamma\neq \emptyset$. So $F_{k}g\gamma^{-1}\cap F_{m}\neq\emptyset$. By the induction hypothesis, we have $F_{k}g\gamma^{-1}\cap F_{k}D_{k}^{m}\neq \emptyset$ and hence $F_{k}g\cap F_{k}D_{k}^{m}\gamma\neq \emptyset$. By the definition of $D_{k}^{m}$, we have $D_{k}^{m}\gamma\subset D_{k}^{m}\delta_{m}^{n}\subset D_{k}^{n}$. Thus $F_{k}g\cap F_{k}D_{k}^{n}\neq \emptyset$ and the claim follows from induction.

\medskip
\noindent{\bf Claim 4}. The $D_{k}^{n}$-translates of $F_{k}$ are maximally disjoint within $T_{k}^{n}$ for all $n,k\in\mathbb{N}$ with $k<n$.

\medskip
Fix $k,n\in\mathbb{N}$ with $k<n$ and let $g\in G$ be such that $F_{k}g\subset T_{k}^{n}$.  We need to show $F_{k}g\cap F_{k}D_{k}^{n}\neq\emptyset$. If $F_{k}g\cap F_{k}\delta_{k}^{n}\neq\emptyset$ then we are done, since $\delta_{k}^{n}=D_{k}^{k}\delta_{k}^{n}\subset D_{k}^{n}$. Now we assume that $F_{k}g\cap F_{k}\delta_{k}^{n}=\emptyset$. According to the definition of $\delta_{k}^{n}$, we cannot have $F_{k}g\subset T_{k}^{n}-\bigcup_{k<m<n}F_{m}\delta_{m}^{n}$. Thus  $F_{k}g\cap (\bigcup_{k<m<n}F_{m}\delta_{m}^{n})\neq\emptyset$. Let $m$ with $k<m<n$ and $\gamma\in\delta_{m}^{n}$ be such that $F_{k}g\cap F_{k}\gamma\neq\emptyset$. Then $F_{k}g\gamma^{-1}\cap F_{k}\neq\emptyset$. By Claim 3, we have $F_{k}g\gamma^{-1}\cap F_{k}D_{k}^{m}\neq\emptyset$ and hence $F_{k}g\cap F_{k}D_{k}^{m}\gamma\neq\emptyset$. Since $D_{k}^{m}\gamma\subset D_{k}^{m}\delta_{m}^{n}\subset D_{k}^{n}$, we have $F_{k}g\cap F_{k}D_{k}^{n}$. This shows Claim 4.

\bigskip
By Claim 4, the $D_{k}^{n}$-translates of $F_k$ are maximally disjoint within $H_{n-1}\cap T$ since $H_{n-1}\cap T\subset T_{k}^{n}$.  This establishes (d). Recall that $G=\bigcup_{n\in\mathbb{N}}H_n$. It now follows that $\Delta_k$-translates of $F_k$ are maximally disjoint within $T$.

\end{proof}

Finally we are able to give
\begin{proof}[Proof of Theorem \ref{existence of blueprint}]
Now Theorem \ref{existence of blueprint} follows from Lemma \ref{existence of growth sequence} and Theorem \ref{existence of blueprint under growth sequence}.
\end{proof}

\appendix


\section{Proximal quasi-factor}
In this section, we show the other direction of Lemma \ref{incontractible} not only for the completeness but also for the sake of safety, since there is a bit difference between our settings and that in \cite[Theorem]{GLA}. However, the proof given here is almost the same as in \cite{GLA}.

\medskip
In this section, a dynamical system $(X, T)$ is a compact metric space $X$ endowed with a continuous action of a topological group $T$.

\medskip
First, we recall the definition of $P$-system given in \cite[Section 3]{GLA}. If $(\mathcal{X}, T)$ is a    $T$-system and $\mathcal{X}$ is equipped with a partial ordering $\geq$ satisfying
\begin{itemize}
\item[(a)] If $X_{\nu}$ and $Y_{\nu}$ are convergent nets in $\mathcal{X}$ with $\lim X_{\nu}=X$ and $\lim Y_{\nu}=Y$, and if for each $\nu$, $X_{\nu}\geq Y_\nu$, then $X\geq Y$,
\item[(b)] For any $X, Y\in \mathcal{X}$, if $X\geq Y$, then for every $t\in T$, $tX\geq tY$,
\end{itemize}
then we call $(\mathcal{X}, T)$ is an {\it ordered system}. A subset $\mathcal{Q}$ of an ordered system $(\mathcal{X}, T)$ is called {\it hereditary} if for any $X\in \mathcal{Q}$ we have $\{Y\in \mathcal{X}: Y\leq X\}\subset \mathcal{Q}$. An ordered system $(\mathcal{X}, T)$ is called {\it irreducible} if there is no proper closed invariant hereditary subset of $\mathcal{X}$.

\medskip
Clearly, for any system $(X,T)$, $T$ naturally acts on the hyperspace ${\rm CL}(X)$ consisting of nonempty closed subsets of $X$. Equipping the ordering of inclusion on ${\rm CL}(X)$  (i.e., $A\leq B\Leftrightarrow A\subset B$), the system $({\rm CL}(X), T )$ becomes an ordered system. A minimal subsystem of $({\rm CL}(X), T )$ is called a {\it quasi-factor} of $(X,T)$.

\medskip
Finally, an ordered system $(\mathcal{X}, T)$ is called a {\it $P$-system} if any two elements $X,Y$ of $\mathcal{X}$ have a $T$ {\it lower bound}, that is if there is a net $t_{\nu}$ in $T$ and $X',Y',Z\in\mathcal{X}$ such that
\[\lim t_{\nu}X=X',\quad \lim t_{\nu}Y=Y'\quad \text{ and }\quad X'\geq Z, \quad Y'\geq Z.\]
The next proposition implies that every $P$-system contains a proximal subsystem.
\begin{prop}\cite[Theorem 3.1]{GLA}\label{irreducible subsystem}
If $(\mathcal{X}, T)$ is a $P$-system then there is a unique irreducible subsystem $(\mathcal{Q},T)$ which is a proximal system and contains a unique minimal set.
\end{prop}

Now let $(X,T)$ be a system and $E(X)$ be the enveloping semigroup of $(X,T)$, which is defined to be the closure of $T$  in $X^{X}$ endowed with the pointwise convergence topology. Let $I$ be a minimal left ideal of $E(X)$ and let $J$ be the set of idempotents in $I$. It is known that $I$ is closed in $E(X)$ and $J$ is nonempty (see \cite[Chapter 6]{Aus}). The following properties is needed.
\begin{lem}\label{min ideal}
Let $I$ be a minimal left ideal of the enveloping semigroup $E(X)$ and $J$ be the set of idempotents in $I$. Then
\begin{itemize}
\item[(1)] For any $u\in J$ and $p\in I$, we have $pu=p$;
\item[(2)] For each $u\in J$, $uI$ is a group;
\item[(3)] The set $\{vI: v\in J\}$ forms a partition of $I$;
\item[(4)] If $ux=x$ for some $u\in J$ then $x$ is a minimal point.
\end{itemize}
\end{lem}

For a system $(X,T)$, we define
\[\mathcal{X}=\{ A\in {\rm CL}(X): \exists v\in J \text{ s.t. } vX\subset A\},\]
where $vX=\{vx: x\in X\}$.  The following proposition is shown in \cite{GLA} for minimal systems but actually it does not use the minimality of $(X,T)$.  For completeness, we repeat the proof here.
\begin{prop}\cite[Proposition 3.2]{GLA}\label{P-system}
$(\mathcal{X}, T)$ is a $P$-system.
\end{prop}
\begin{proof} 
First we show $\mathcal{X}$ is closed in ${\rm CL}(X)$. If $A_i$ is a sequence in $\mathcal{X}$ which converges to $A\in {\rm CL}(X)$, then for each $i$, there is a $v_i\in J$ with $v_iX\subset A_i$. We may assume that $v_i\rightarrow p\in I$. By Lemma \ref{min ideal}, there are $v\in J$ and $\alpha\in uI$ such that $p=v\alpha$ and hence $v=p\alpha^{-1}$, since $uI$ is a group. Now for any $x\in X$, $\lim v_i(\alpha^{-1}x)=p\alpha^{-1}x=vx$. Thus $vX\subset A$ and hence $A\in \mathcal{X}$. This shows that $\mathcal{X}$ is closed in ${\rm CL}(X)$.

\medskip
Now we show that $\mathcal{X}$ is invariant under $T$ action.  For this, fix a $t\in T$ and $A\in \mathcal{X}$. Then $A\supset vX$ for some $v\in J$. Similarly, it follows from Lemma \ref{min ideal} that $tv=w\alpha$ for some $w\in J$ and $\alpha\in uI$.  Then for each $x\in X$, $wx=(w\alpha)(\alpha^{-1}x)=(tv)(\alpha^{-1}x)$ so that $tA\supset tvX=w\alpha X=wX$. So $\mathcal{X}$ is invariant.

\medskip
Finally, we show $(\mathcal{X},T)$  is a $P$-system. $(\mathcal{X},T)$ is an ordered system as a subsystem of $({\rm CL}(X), T)$. It remains to show any two elements of $\mathcal{X}$ have $T$ lower bound. Let $A,B$ be in $\mathcal{X}$ then there are $v, w\in J$ with $vX\subset A$ and $wx\subset B$. Let $t_{i}$ be a net in $T$  with
\[\lim t_i=u,\quad \lim t_i A=A',\quad \lim t_i B=B'.\]
Then we have $uX\subset A'\cap B'$ which implies that $A$ and $B$ have a $T$ lower bound in $\mathcal{X}$. Therefore, $(\mathcal{X},T)$  is a $P$-system.
\end{proof}

Now it follows from Proposition \ref{irreducible subsystem} and Proposition \ref{P-system} that there is a unique proximal quasi-factor $(\mathcal{X}_0, T)$ in $(\mathcal{X}, T)$, which we denote by $(\Pi(X), T)$. Actually, the definition of $\mathcal{X}$ and hence $\Pi(X)$ does not depend on the choice of the minimal left ideal $I$ in $E(X)$ and the idempotent $u\in I$ (see \cite[section 3]{GLA}). For our aim, we just need to fix such an $I$ and $u$.

\begin{defn}
We say a system $(X,T)$ is {\it non-contractible} if $(\Pi(X), T)$ is the trivial system, i.e., $\Pi(X)=\{X\}$.
\end{defn}
\begin{lem}\label{Pi(X)}
If the set of minimal points in $(X,T)$ is dense, then for any idempotent $u\in J$, we have
\[ X\subset \bigcup\left\{ A: A\in \overline{T(\overline{uX}^{X})}^{{\rm CL}(X)}\right\},\]
where $\overline{ \ * \ }^X$ and $\overline{ \ * \ }^{{\rm CL}(X)}$ stand for the closure taken in $X$ and ${\rm CL}(X)$, respectively.
\end{lem}
\begin{proof}
Let $M$ denote the set of minimal points of $(X,T)$. Then for any $x\in X$, there is a sequence $(x_n)$ in $M$ converging to $x$. Now for  each $n$, $ux_n\in \overline{Tx_{n}}^{X}$ that is also a minimal point. Thus, for each $n$, there is a $t_n\in T$ such $\rho(t_nux_n, x_n)<\frac{1}{2^{n}}$.
It follows that $\lim t_nux_n=\lim x_n=x$. By the definition of the hyperspace topology, we conclude the proof.
\end{proof}

Recall that a system $(X,T)$ is incontractible if the set of minimal points in $X^{n}$  is dense for each $n\geq 1$.
\begin{lem}\label{non>in}
If $(X,T)$ is non-contractible then it is incontractible.
\end{lem}
\begin{proof}
Let $I$ be the minimal left ideal in $E(X)$ to define $\mathcal{X}$ and $J$ be the set of idempotents of $I$. If $\Pi(X)=\{X\}$, then for each $v\in J$ there is a net $t_i$ in $T$ such that  $\lim t_i(\overline{vX})=X$ in ${\rm CL}(X)$.  Indeed, otherwise there is some $v\in J$ such that $\mathcal{X}_1=\{t(\overline{vX}): t\in T\}$ is a subsystem in $\mathcal{X}$ which does not contain $X$; and then we can find a proximal quasi-factor in $\mathcal{X}_1$. This contradicts our assumption that $(X,T)$ is non-contractible.

\medskip
Now fix an $u\in J$ and $n\geq 1$. Then for any $x_1,x_2,\cdots, x_n\in X$, we have
\[ u(ux_1, \ldots, ux_n) =(ux_1,\ldots, ux_n).\]
By Lemma \ref{min ideal}, $(ux_1,\ldots, ux_n)$ is a minimal point in $X^{n}$. From the above paragraph,  there is a net $t_i$ in $T$ such that $\lim t_i(\overline{uX})=X$ in ${\rm CL}(X)$.
Thus for any $(y_1,\ldots, y_n)\in X^{n}$, there are $z_1,\ldots,z_n\in \overline{uX}$ such that  $\lim (t_iz_1,\ldots, t_iz_n)=(y_1,\cdots, y_n)$. Further, there are $x_1,\ldots, x_n\in X$ such that $\lim (t_iux_1,\ldots, t_iux_n)=(y_1,\cdots, y_n)$. Note that for each $i$, $(t_iux_1,\ldots, t_iux_n)$ is a minimal points in $X^{n}$. Thus $(X,T)$ is incontractible.
\end{proof}

Now we are ready to prove the other direction of  Lemma \ref{incontractible} as done in \cite[Theorem 4.2]{GLA}.
\begin{lem}
If $(X,T)$ has a dense set of minimal point and is disjoint from every minimal proximal $T$-system then it is incontractible.
\end{lem}
\begin{proof}
Let $(\Pi(X), T)$ be the proximal quasi-factor of $(X,T)$. Then $(X,T)$ is disjoint with $(\Pi(X), T)$. Let
\[J=\{(x, A): A\in \Pi(x), x\in A\}\subset X\times \Pi(X).\]
Clearly,  the projection of $J$ to $\Pi(X)$ is full. By Lemma \ref{Pi(X)},  the projection of $J$ to $X$ is also full since $\Pi(X)$ is an invariant closed subset in ${\rm CL}(X)$. Thus $J$ is a joining of $X$ and $\Pi(X)$. By the disjointness, we have $J=X\times \Pi(X)$ and hence $\Pi(X)=\{X\}$. Now it follows from Lemma \ref{non>in} that $X$ is incontractible.
\end{proof}

\bigskip
\noindent{\bf Acknowledgement}:  We are grateful for the careful reading and useful comments of the anonymous referee. We also thank Su Gao, Wen Huang and Song Shao for their useful discussions.


\begin{thebibliography}{HD}

\bibitem{AG} E. Akin and E. Glasner, {\it Residual properties and almost equicontinuity}. J. Anal. Math. {\bf 84} (2001), 243-286.

\bibitem{AK} E. Akin and S. Kolyada, {\it Li-York sensitivity}. Nonlinearity {\bf 16} (2003), 1421-1433.

\bibitem{Aus} J. Auslander, {\it Minimal flows and their extensions}. North-Holland Mathematics Studies {\bf 153} (North-Holland, Amsterdam), 1988.

\bibitem{BF} B. Balcar and F Franek, {\it Structureal prperties of universal minimal dynamical systems for discrete seigroups}. Trans. Amer. Math. Soc. {\bf 349}(5) (1997), 1697-1724.

\bibitem{Ber} A. Bernshteyn, {\it A short proof of Bernoulli disjointness via local lemma}. Proc. Amer. Math. Soc. {\bf 148} (2020), no. 12, 5235-5240.

\bibitem{Bollobas} B. Bollob\'{a}s, {\it Moder graph theory}. GTM {\bf 184}, Springer, 1998.

\bibitem{CKN} G. Cairns, A. Kolganova and A. Nielsen, {\it Topological transitivity and mixing notions for group actions}. Rocky Mountain J. Math. {\bf 37} (2007), no. 2, 371-397.


\bibitem{CP} M. I. Cortez and S. Petite, {\it $G$-odometers and their almost one-to-one extensions}. J. London Math. Soc. (2) {\bf 78} (2008), 1-20.

\bibitem{DSY} P. Dong, S. Shao and X. Ye, {\it Product recurrent properties, disjointness and weak disjointness}. Israel J. Math. {\bf 188} (2012), 463-507.

\bibitem{Ellis} R. Ellis, {\it Universal minimal sets},  Proc. Amer. Math. Soc. {\bf 11} (1960), 540-543.

 \bibitem{Fur} H. Furstenberg,  {\it Disjointness in ergodic theory, minimal sets, and a problem in Diophantine approximation}, Math. Systems Theory {\bf 1} (1967), 1-49.

 \bibitem{GJS} S. Gao, S. Jackson and B. Seward, {\it Group colorings and Bernoulli subflows}. Mem. Amer. Math. Soc. {\bf 241}(1141), May 2016.

\bibitem{GLA} S. Glasner, {\it Compressibility properties in topological dynamics}. Amer. J. Math. {\bf 97} (1975), 148-171.

 \bibitem{GTWZ} E. Glasner, T. Tsankov, B. Weiss and A. Zucker, {\it Bernoulli disjointness}. Duke Math. J. {\bf 70} (2021) no. 4,  614-651.

\bibitem{HSY} W. Huang, S. Shao and X. Ye, {\it An answer to Furstenberg's problem on topological disjointness}. Ergod. Th. \& Dynam. Sys. {\bf 40} (2020), 2467-2481.

\bibitem{HY} W. Huang and X. Ye, {\it Dynamical systems disjoint from any minimal system}. Trans. Amer. Math. Soc. {\bf 357} (2) (2005), 669-694.

\bibitem{KRS} M. Kennedy, S. Raum and G. Salomon, {\it Amenability, proximality and higher-order syndeticity}. Forum of Mathematics, Sigma {\bf Vol. 10: e22} (2022), 1-28.



\bibitem{LYY}J. Li, K. Yan and X. Ye,  {\it Recurrence properties and disjointness on the induced spaces}. Discrete Contin. Dyn. Syst.
{\bf 35}(2015),  1059-1073.

\bibitem{MZ} D. Montgomery and L. Zippin, {\it Topological transformation groups}. Interscience Publishers, New York-London, 1955. xi+282 pp.

\bibitem{Opr10} P. Oprocha, {\it Weak mixing and product recurrence}. Ann. Inst. Fourier (Grenoble) {\bf 60} (4) (2010), 1233-1257.

\bibitem{Opr19} P. Oprocha, {\it Double minimality, entropy and disjointness with all minimal systems}. Discrete Contin. Dyn. Syst. {\bf 39} (2019), 263-275.

\bibitem{Roman}S.  Roman,  {\it Fundamentals of group theory. An advanced approach} . Birkh\"{a}user/Springer, New York, 2012. xii+380 pp.

\bibitem{Veech0} W. Veech, {\it Point-distal flows}. Amer. J. Math. {\bf 92} (1970), 205-242.

\bibitem{Veech} W. Veech, {\it Topological dynamics}, Bull. Amer. Math. Soc. {\bf 83}(5) (1977), 775-830.



\bibitem{WCF} H. Wang, Z. Chen and H. Fu, {\it $M$-systems and scattering systems of semigroup actions}. Semigroup Forum {\bf 91} (2015), 699-717.

\bibitem{Yu} T. Yu, {\it Dynamical systems disjoint from any minimal system under group actions}. Difference Equations, Discrete Dynamical Systems and Applications, 181-195 (Springer Proceedings in Mathematics and Statistics, {\bf 150}). Springer, Cham, 2015.



\end{thebibliography}
\end{document}